\documentclass[journal]{IEEEtran}
\usepackage{amsmath,amsfonts,amsthm}
\usepackage{array}
\usepackage[caption=false,font=normalsize,labelfont=sf,textfont=sf]{subfig}
\usepackage{textcomp}
\usepackage{stfloats}
\usepackage{url}
\usepackage{microtype} 
\usepackage[subtle,tracking=normal]{savetrees} 
\usepackage{hyperref} 
\usepackage{verbatim}
\usepackage{graphicx}
\usepackage{xcolor}
\usepackage{nomencl}
\usepackage{etoolbox}
\usepackage{ifthen}
\usepackage{booktabs}
\usepackage[ruled]{algorithm2e}
\usepackage{placeins}
\usepackage{multirow}
\usepackage{makecell}
\usepackage{enumitem}
\usepackage{pdfpages}

\allowdisplaybreaks 

\makenomenclature
\usepackage{cite}
\hypersetup{
    colorlinks=true, 
    linkcolor=blue, 
    filecolor=black, 
    urlcolor=blue, 
    citecolor=blue, 
    linkbordercolor={1 1 1}, 
    pdfborder={0 0 0} 
}

\def\BibTeX{{\rm B\kern-.05em{\sc i\kern-.025em b}\kern-.08em
    T\kern-.1667em\lower.7ex\hbox{E}\kern-.125emX}}
\usepackage{balance}

\newtheorem*{remark}{Remark}


\begin{document}
\bstctlcite{IEEEexample:BSTcontrol}

\newpage
\title{
Quantum Learning and Estimation for \\Coordinated Operation between Distribution Networks\\ and Energy Communities 
}

\author{
\IEEEauthorblockN{Yingrui Zhuang}, 
\IEEEauthorblockN{Lin Cheng, \textit{Senior Member, IEEE}},
\IEEEauthorblockN{Yuji Cao},
\IEEEauthorblockN{Tongxin Li, \textit{Member, IEEE}}, \\
\IEEEauthorblockN{Ning Qi, \textit{Member, IEEE}},
\IEEEauthorblockN{Yan Xu, \textit{Senior Member, IEEE}},
\IEEEauthorblockN{Yue Chen, \textit{Senior Member, IEEE}}

\thanks{
This work was supported in part by the 1+1+1 CUHK-CUHK(SZ)-GDST Joint Collaboration Fund (No. 2025A0505000047), 
National Natural Science Foundation of China (No. 52037006), the CUHK Strategic Seed Funding for Collaborative Research Scheme (No. 3136080),
the CUHK Strategic Partnership Award for Research Collaboration (No. 4750467), 
and the China Postdoctoral Science Foundation special funded project (No. 2023TQ0169). (\textit{Corresponding author}: Yue Chen.)

Yingrui Zhuang and Lin Cheng are with the Department of Electrical Engineering, 
Tsinghua University, Beijing 100084, China (email: zyr21@mails.tsinghua.edu.cn).

Yuji Cao and Yue Chen are with the Department of Mechanical and Automation
Engineering, The Chinese University of Hong Kong, Hong Kong, China (email: yjcao@mae.cuhk.edu.hk, yuechen@mae.cuhk.edu.hk).

Tongxin Li is with the School of Data Science, The Chinese University of Hong Kong, Shenzhen, China (email: litongxin@cuhk.edu.cn).

Ning Qi is with the Department of Earth and Environmental Engineering, Columbia University, NY 10027, USA (email: nq2176@columbia.edu). 

Yan Xu is with the School of Electrical and Electronic Engineering, Nanyang Technological University, Singapore 639798, Singapore (email: xuyan@ntu.edu.sg).}
}

\markboth{CSEE JOURNAL OF POWER AND ENERGY SYSTEMS, VOL. XX, NO. XX, XX 2025}
{How to Use the IEEEtran \LaTeX \ Templates}

\maketitle

\begin{abstract}
    Price signals from distribution networks (DNs) guide energy communities (ECs) in adjusting their energy usage, enabling effective coordination for reliable power system operation.
    However, this coordinated operation faces significant challenges due to the limited availability of ECs' internal information (i.e., only the aggregated energy usage of ECs is available to DNs),
    and the high computational burden of accounting for uncertainties and the associated risks through numerous scenarios.   
    To address these challenges, we propose a quantum learning and estimation approach to enhance coordinated operation between DNs and ECs.
    Specifically, by
    leveraging advanced quantum properties such as quantum superposition and entanglement,
    we develop a hybrid quantum temporal convolutional network-long short-term memory (Q-TCN-LSTM) model
    to establish an end-to-end mapping between ECs' responses and the price incentives from DNs. 
    Moreover, we develop a quantum estimation method based on quantum amplitude estimation (QAE) and two phase-rotation circuits
    to significantly accelerate the optimization process under numerous uncertainty scenarios.
    Numerical experiments demonstrate that,
    compared to classical neural networks,
    the proposed Q-TCN-LSTM model improves the mapping accuracy by 69.2\% while reducing the model size by 99.75\%.
    Compared to classical Monte Carlo simulation,
    QAE achieves comparable accuracy with a substantial reduction in computational resources. 
    In addition, the estimated computation time for quantum learning and estimation on ideal quantum devices is over 90\% shorter than that of traditional methods.
\end{abstract}

\mbox{}

\vspace{-6mm}
\section{Introduction}\label{sec:introduction}

\IEEEPARstart{W}{ith} the rapid development of distributed energy resources (DERs),
traditional centralized power systems are gradually evolving toward a more decentralized architecture~\cite{xu2024coordinated,xie2024collaborative}.
Energy communities (ECs), which integrate various types of DERs, facilitate local consumption of renewable energy and enhance the local energy balance~\cite{Networked_EC_review}.
However, from the standpoint of upper-level distribution networks (DNs), the operational behaviors of ECs inherently exhibit uncertainties.
In this context, achieving coordinated operation between DNs and ECs 
is essential for enhancing system flexibility and reducing potential risks (e.g., voltage issues, congestion)~\cite{DNECcoordination4}.
Among existing coordination methods, one effective and commonly used approach is to leverage price signals from DNs to guide grid-supportive responses from ECs~\cite{incentative}.
However, two critical challenges remain in achieving effective DNs-ECs coordinated operation.

The first challenge lies in accurately estimating ECs' electricity consumption behavior~\cite{Consumption_Pattern,zhou2024forecasting}. 
In the existing literature, an ideal assumption is often adopted that DNs possess full access to the internal operational data of ECs, including detailed device parameters and operational preferences. 
Consequently, precise operational models for ECs can be developed based on the physical characteristics of individual devices, such as the electro-thermal coupling dynamics of air conditioners~\cite{AC} and specific industrial production processes~\cite{Industrial3}. 
Moreover, to improve the overall computational tractability,
simplified convex EC operation models are built 
and embedded into the DNs operation model to form
a centralized coordinated operation model~\cite{SG1}.
However, this ideal assumption rarely holds true in practice for two primary reasons.
First, 
user electricity consumption behaviors usually exhibit strong nonlinearity, 
non-stationarity, and partial irrationality, 
which cannot be accurately captured by simplified convex operational models. 
Second, ECs are unwilling or unable to share detailed, device-level data due to concerns about user privacy, cybersecurity risks, or commercial confidentiality~\cite{privacy2}, making it impractical to develop precise operational models based on internal data. 
Consequently, in practical applications, distribution networks usually have access only to aggregated energy consumption data from ECs without detailed insights into device-level configurations or usage patterns~\cite{smartmeter}. 
To address this issue, 
distributed optimization methods (e.g., alternating direction method of multipliers~\cite{ADMM}) have been widely adopted, 
where the centralized optimization problem is decomposed into several subproblems that are solved independently. While such methods offer good privacy protection, 
they often suffer from heavy communication burdens and slow convergence. Alternatively,
some studies have proposed data-driven approaches to construct surrogate models that emulate user electricity consumption behaviors in response to external incentives~\cite{ECfit}. 
However, such behavioral modeling~\cite{qi2023chance} is essentially a time-series analysis problem characterized by high dimensionality, temporal coupling, and
strong nonlinearity, leading to high learning complexity. 
In particular, 
for ECs with complex operational characteristics, 
achieving high behavioral modeling accuracy often requires a large number of model parameters~\cite{neural_scaling_laws}, 
which further increases model complexity and reduces computational efficiency.

The second challenge arises from the substantial computational burden imposed by risk-aware models that rely heavily 
on a large number of discrete scenarios to accurately characterize uncertainties and potential operational risks, 
especially when precise probability distributions are unavailable~\cite{SAA}.
This issue becomes especially pronounced when addressing extreme loss events using risk measures such as conditional value-at-risk (CVaR), which require accurate representations of distribution tails, thus demanding extensive scenario sampling to capture rare but high-impact events~\cite{cvar5}. Typically, Monte Carlo (MC) simulation is employed to evaluate the risk expectation across a broad range of possible scenarios~\cite{MC}. 
Nevertheless, 
MC’s estimation accuracy is correlated with the number of sample scenarios
and often suffers from low computational efficiency and slow convergence. Scenario reduction techniques have been explored to compress the original large-scale scenario set into a smaller yet
representative subset~\cite{SR} to reduce computational complexity. 
However, 
for tail-sensitive risk measures such as CVaR, 
scenario reduction may not be able to fully preserve the integrity of distribution tails, 
leading to a loss of critical information and degraded risk-awareness.

Recent advances in quantum computing have opened up new opportunities for effective DNs-ECs coordinated operation.
Unlike classical computing, which processes information sequentially and requires large parameter sets to model complex systems,
quantum computing leverages quantum superposition and entanglement properties inherent in quantum bits (qubits),
allowing the representation and manipulation of an exponential number of states simultaneously with fewer parameters.
On the hardware front, companies such as IBM and Google have developed prototype systems with tens to hundreds of qubits, 
signaling the advent of the Noisy Intermediate-Scale Quantum (NISQ) era~\cite{quantum_advances1}. 
These developments make quantum computing particularly suitable for solving problems characterized by high dimensionality, complex data structures, and intensive computational demands with inherent parallelism.

To address the first challenge regarding the accurate estimation of ECs' electricity consumption behavior, 
quantum neural networks (QNNs) present a promising approach for modeling the mapping between external incentives and EC responses. Compared to classical neural networks~\cite{chen2023graph}, 
QNNs leverage qubits and quantum gates to perform computations within superposed and entangled quantum states, 
potentially enhancing their capability to represent complex nonlinear relationships and accelerating the training process. 
Existing studies have demonstrated the effectiveness of QNNs in various power system applications, 
including carbon price forecasting~\cite{Qforecast1} and fault diagnosis~\cite{QDL1}. 
However, the incentive-response modeling of ECs involves heterogeneous operational characteristics, 
forming a significantly more intricate temporal regression task. 
The QNN architectures employed in current research remain relatively simplified, 
limiting their ability to fully capture the complex temporal dynamics and coupled operational dependencies inherent in EC behaviors.

For the second challenge of computational bottlenecks arising from large-scale scenario evaluation,
quantum estimation methods (e.g., quantum amplitude estimation, QAE) provide a promising solution.
Compared to classical MC methods, QAE theoretically achieves a quadratic speedup~\cite{MC1} under specific conditions (e.g., efficient quantum state preparation and negligible quantum noise), 
demonstrating significant advantages in key subroutines such as expectation evaluation and CVaR computation~\cite{QAE1}. 
Refs~\cite{QAEreliability,QAEreliability1} have applied QAE to power system reliability assessment, 
validating its potential in uncertainty-aware analysis. However, existing research predominantly addresses theoretical aspects of QAE algorithms, with limited exploration into practical quantum circuit implementations. Consequently, a significant challenge is designing quantum circuits capable of effectively balancing estimation accuracy and computational complexity, thus enabling rapid and precise uncertainty scenario evaluation in practical applications.

In this paper, we propose a quantum learning and estimation method for the coordinated operation of DNs and ECs.
Our main contributions are as follows:
\begin{enumerate}
  \item We develop a quantum temporal convolutional network–long short-term memory (Q-TCN-LSTM) model  
        to estimate the energy consumption behavior of ECs in response to the price incentives from the DNs. 
        To the best of our knowledge, the Q-TCN-LSTM model has not been explored in prior research.
        Q-TCN and Q-LSTM respectively extract short-term and long-term temporal dependencies,
        comprehensively modeling the coupled temporal dependencies of incentive-response characteristics of ECs.
        Compared to the classical TCN-LSTM,
        Q-TCN-LSTM improves accuracy by 69.2\%,
        while reducing the model size by 99.75\%.
        Moreover, on ideal quantum devices, the theoretically estimated learning time for one forward pass is over 90\% shorter than that required by traditional learning methods.
    \item We develop a quantum estimation method to address the computational burden of estimating the gradient of CVaR-based optimization under a large number of uncertainty scenarios, which is a critical component in DNs–ECs coordinated operation.
    Two phase rotation circuits are introduced: 
    one establishes analytical linear approximations between rotation angles and quantum states, 
    achieving low circuit complexity but with potential approximation errors; 
    the other computes exact rotation angles which are adaptively applied by introducing additional qubits, 
    ensuring high accuracy but with increased complexity. 
    For the proposed quantum estimation method, the theoretically estimated computation time is reduced by up to 99.99\% compared with classical Monte Carlo simulation.

\end{enumerate}

The rest of the paper is organized as follows. 
Section~\ref{sec:prob} introduces the coordinated operation problem formulation between DNs and ECs.
The methodology of the proposed quantum learning and estimation approach is presented in Section~\ref{sec:methodology}.
Numerical case studies are provided in Section~\ref{sec:case}. 
Finally, we summarize our conclusions in Section~\ref{sec:conclusion}.

\section{Problem Formulation}\label{sec:prob}
In this section, we first present the DNs-ECs coordinated operation structure, 
followed by the detailed risk-averse operation model of DNs,
and the coordination mechanism between DNs-ECs coordinated operation.

\subsection{DNs-ECs Coordinated Operation Structure}
\begin{figure}[htbp] 
  \setlength{\abovecaptionskip}{-0.1cm}  
  \setlength{\belowcaptionskip}{-0.1cm}   
  \centering
  \includegraphics[width=1.00\columnwidth]{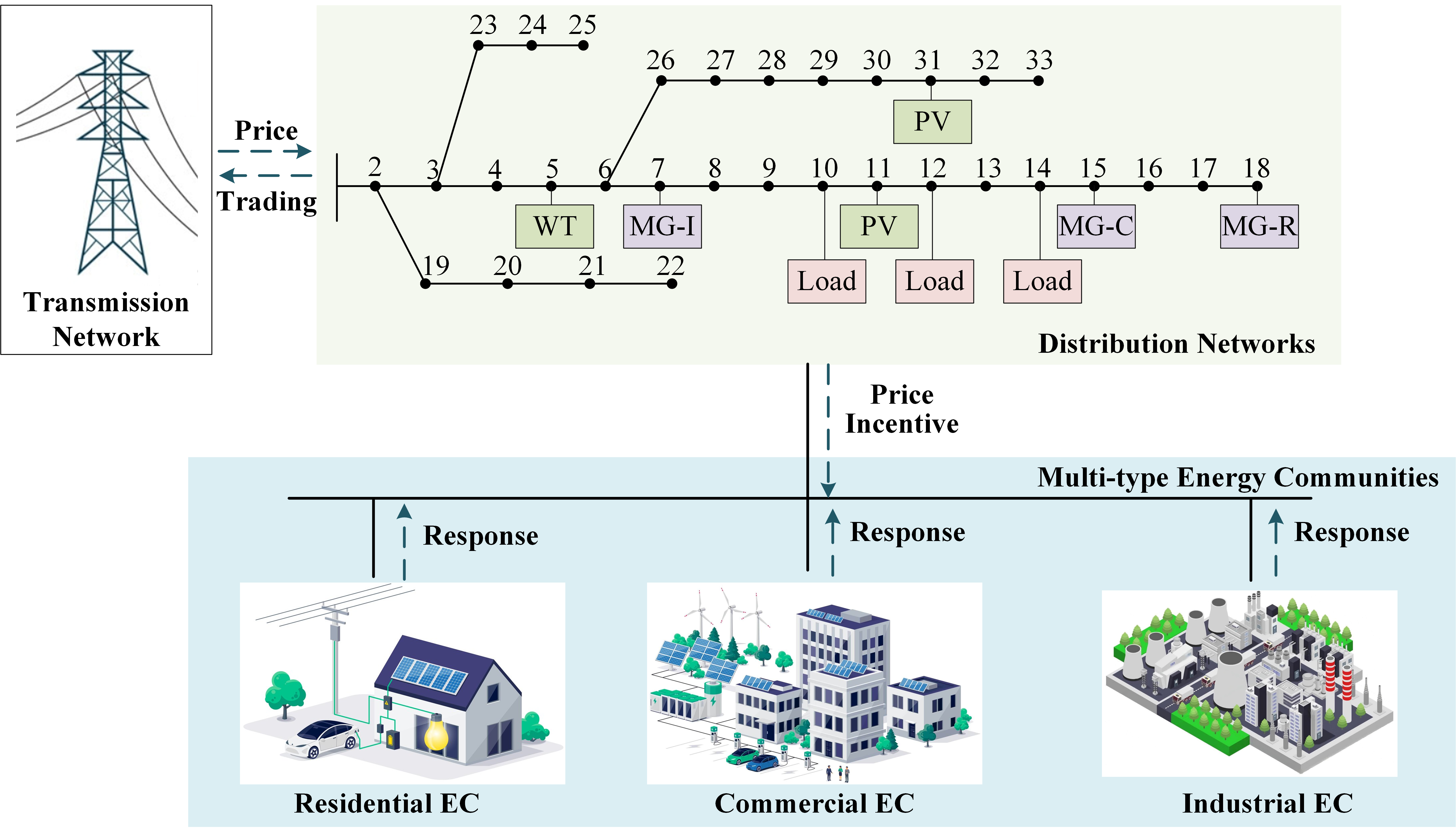}
  \caption{Coordinated operation structure between DNs and multiple ECs.}
  \label{fig:coordination}
\end{figure}

The structure of DNs integrated with multiple ECs 
and DERs (e.g., wind turbines (WTs), photovoltaic (PV) systems, and new-type loads)
is illustrated in Fig.~\ref{fig:coordination}. 
This paper focuses on three representative types of ECs: residential (EC-R), commercial (EC-C), and industrial (EC-I).
These three types not only cover the main categories of ECs in real-world situations thus ensuring the generalizability,
but also exhibit significant heterogeneity across multiple dimensions
(e.g., device composition, operational constraints and objectives, risk preferences), leading to distinct operational behaviors and introducing
non-trivial complexity into DNs-ECs coordinated operation.

In this paper, we focus on the day-ahead price-incentive demand response mechanism as the coordinated operation strategy,
where price signals guide ECs to adjust their operations accordingly.
The coordinated operation between DNs and ECs is naturally a bilevel optimization problem,
with DNs making decisions at the upper level and ECs responding at the lower level.
At the upper level, DNs should accurately estimate the uncertainty associated with renewable energy sources (RESs), loads, and the aggregated energy consumption behaviors (i.e., trading power) of ECs,
and accordingly determine risk-averse operational decisions, which is mainly considered as the price signals for ECs.
At the lower level, 
each EC independently optimizes its operational strategy based on customized operational objectives
in response to the price signals from the DNs,
and submits the resulting power trading plans to the DNs.

\subsection{Risk-Averse Operation Model of Distribution Networks}
Considering the uncertainties in load and RESs, 
DNs need to make risk-averse decisions to reduce the operational risks
including voltage violations and line overloads.
In this paper, CVaR is adopted to quantify operational risks, where the uncertainties in load and RESs are modeled using a set of discrete scenarios.
Accordingly, the overall operational model of the DNs is formulated as follows:
\begin{subequations}\label{DNs_opt}
  \begin{align}
    \label{obj_all}
    \min_{\pi^{\mathrm{EC}},v_{\alpha}}\ &C= \sum_{s=1}^{S}\gamma_sC_s + \lambda ( v_{\alpha} + \frac{1}{1-\alpha} \sum_{s=1}^{S}\gamma_s  [ C_s -v_{\alpha} ]_{+} ), \\
    \text{s.t.} \ 
    \label{objs}
    & C_s = \sum_{t=1}^{T} \Delta t[ -\sum_{i\in \Omega^{\mathrm{EC}}} \pi^{\mathrm{EC}}_{t}  P^{\mathrm{EX,EC}}_{i,t}  + \pi^{\mathrm{DN}}_{t}  P^{\mathrm{EX,DN}}_{s,t}  \nonumber \\
    & \quad + \beta_1    \sum_{i\in \Omega^{\mathrm{B}}}([V_{i,s,t} - \overline{V}_{i}]_+ + [\underline{V}_{i}-V_{i,s,t}]_+)   \\
    & \quad + \beta_2  \sum_{ij\in \Omega^{\mathrm{L}}} ([P_{ij,s,t} - \overline{P}_{ij}])_++ [\underline{P}_{ij}-P_{ij,s,t}]_+)],\nonumber\\
    \label{pf U}
    & \mathcal{V}_{j,s,t} = \mathcal{V}_{i,s,t}  - 2(r_{ij}P_{ij,s,t} + x_{ij}Q_{ij,s,t}), \\
    \label{pf P}
    & p_{j,s,t} = P_{ij,s,t}-\sum_{l:j\rightarrow l}P_{jl,s,t} ,\\
    \label{pf Q}
    & q_{j,s,t} = Q_{ij,s,t}-\sum_{l:j\rightarrow l}Q_{jl,s,t}, \\
    \label{pin}
    & p_{j,s,t} = -P^{\mathrm{R}}_{j,s,t} + P^{\mathrm{L}}_{j,s,t} + P^{\mathrm{EX,EC}}_{j,t}, \\
    \label{qin}
    & q_{j,s,t} = Q^{\mathrm{L}}_{j,s,t} + Q^{\mathrm{EX,EC}}_{j,t}, \\
    \label{P_EX_EC}
    & P^{\mathrm{EX,EC}}_{i} = \phi_{i} (\pi^{\mathrm{EC}} ),\\
    \label{price1}
    & \underline{\pi}^{\mathrm{EC}}_{t} \leq \pi^{\mathrm{EC}}_{t} \leq \overline{\pi}^{\mathrm{EC}}_{t},\\
    \label{price2}
    & \frac{1}{T}\sum_{t=1}^{T} \pi^{\mathrm{EC}}_{t} =   \pi^{\mathrm{EC,av}}_{t},
  \end{align}
\end{subequations}
where 
the released incentive price $\pi^{\mathrm{EC}}$ 
and the VaR $v_{\alpha}$ are the decision variables of DNs.~\eqref{pf U}-\eqref{pf Q} are the linearized distflow model.
\eqref{pin}-\eqref{qin} are the power injection constraints.
\eqref{price1}-\eqref{price2} are the price signal constraints.
$C_s$ is the cost under the $s$-th scenario ($S$ in total) of uncertain DERs integrated to the DNs
with probability $\gamma_s$. 
$P^{\mathrm{EX,EC}}_{i}$ is the active trading power between the EC at node $i$ and DNs,
which is the EC's response to the price signal $\pi^{\mathrm{EC}}$.
Note that the DNs do not know the exact operation patterns and device parameters of ECs,
we adopt an unknown function $\phi_{i}(\cdot)$ to represent the corresponding response characteristics, as in~\eqref{P_EX_EC}.
The active and reactive trading power of ECs are assumed to be proportional.
$\lambda$ is the risk aversion coefficient which controls the risk acceptance level.
$[z]_{+} = \max(0,z)$.
$\pi^{\mathrm{DN}}_{t}$ is the trading price published by the transmission network to DNs.
$P^{\mathrm{EX,DN}}_{s,t}$ is the trading power of the DNs with the transmission network under scenario $s$ at time $t$. 
$\pi^{\mathrm{EC}}_{t}$ is bounded by $\overline{\pi}^{\mathrm{EC}}_{t}/\underline{\pi}^{\mathrm{EC}}_{t}$,
and averaged by $\pi^{\mathrm{EC,av}}_{t}$.
$\beta_1/\beta_2$ are the penalty coefficients of the voltage and line power flow violations.
$V_{i,s}$ is the voltage magnitude bounded by $\overline{V}_{i}/\underline{V}_{i}$.
$P_{ij,s,t}/Q_{ij,s,t}$ are the active/reactive power flow of line $ij$.
$r_{ij}/x_{ij}$ are the resistance/reactance.
$p_{j,s,t}/q_{j,s,t}$ are the active/reactive power injection.
$P^{\mathrm{R}}_{j,s,t}$ is the active power output of RESs at bus $j$ under scenario $s$ at time $t$.
$P^{\mathrm{L}}_{j,s,t}/Q^{\mathrm{L}}_{j,s,t}$ are the active/reactive power loads.
$\Omega^{\mathrm{EC}}/\Omega^{\mathrm{B}}/\Omega^{\mathrm{L}}$ denote the sets of nodes with EC, nodes, and lines, respectively.
In~\eqref{pin}, $P^{\mathrm{R}}_{j,s,t}=0$ and $P^{\mathrm{EX,EC}}_{j,s,t}=0$ if bus $j$ is not connected to any RES or EC, respectively.

\subsection{DNs-ECs Coordinated Operation}
Under limited information availability,
DNs can only access the aggregated energy usage of each EC,
while the detailed internal data and operational states of ECs remain unobservable. 
Therefore, the bilevel DNs-ECs coordinated operation problem cannot be solved directly.
To address this challenge, we propose to learn direct mappings that capture the incentive-response behavior of each EC
based on historical data, 
as in~\eqref{eq:map}.
\begin{equation}\label{eq:map}
  \hat{P}^{\mathrm{EX,EC}}_i = \varphi_i(\pi^{\mathrm{EC}};\varTheta_i)\sim \phi_{i} (\pi^{\mathrm{EC}}),
\end{equation}
where the parameters $\varTheta_i$ can be optimized by minimizing a loss function (e.g., mean square error, MSE):
\begin{equation}
  \varTheta_i^* = {\arg\min}_{\varTheta_i}\  \frac{1}{T}\sum_{t=1}^{T} \left(P^{\mathrm{EX,EC}}_{i,t} - \hat{P}^{\mathrm{EX,EC}}_{i,t}\right)^2.
\end{equation}

Then, under the condition that the mapping model accurately captures the incentive-response behavior of ECs (i.e., the loss value is small),
we can utilize the mapping model in~\eqref{eq:map} 
to approximate the unknown response function in~\eqref{P_EX_EC}.
The original bilevel DNs-ECs coordinated operation problem is reformulated as a single-level optimization problem,
yielding solutions approximate to those of the original bilevel formulation.

By exploiting the differentiability of the mapping model~\eqref{eq:map} together with the CVaR-based formulation, in which auxiliary variables and small perturbations are introduced to convexify and eliminate the non-differentiability caused by the $[\cdot]_+$ operator, the integrated distribution network optimization model becomes fully differentiable and can therefore be effectively solved using gradient-based methods.
The gradient descent update rule is:
\begin{subequations}\label{eq:grad}
  \begin{align}
    &\pi^{\mathrm{EC},k+1} = \pi^{\mathrm{EC},k} - \eta \frac{\partial C}{\partial \pi^{\mathrm{EC},k}}, \ \ v_{\alpha}^{k+1} = v_{\alpha}^k - \eta \frac{\partial C}{\partial v_{\alpha}^k},\\
    &\frac{\partial C}{\partial \pi^{\mathrm{EC},k}} = \sum_{s=1}^{S} \gamma_s \frac{\partial C_s}{\partial \pi^{\mathrm{EC},k}} + \frac{\lambda}{1-\alpha} \sum_{s=1}^{S} \gamma_s 
    \frac{\partial C_s}{\partial \pi^{\mathrm{EC},k}}  \mathbb{I}(C_s > v_{\alpha}^{k}),\\
    &\frac{\partial C}{\partial v_{\alpha}^{k}} = \lambda \left( 1 - \frac{1}{1-\alpha} \sum_{s=1}^{S} \gamma_s \mathbb{I}(C_s > v_{\alpha}^{k}) \right),
    \end{align}
\end{subequations}
where $\mathbb{I}(C_s > v_{\alpha})=1$ if $C_s > v_{\alpha}$ and 0 otherwise.
The feasible region of $\pi^{\mathrm{EC}}$ is ensured by repeatedly doing: 
firstly normalized to [0,1] via softmax then inverse-mapped to actual range to satisfy range constraint~\eqref{price1}
and mean-shifted to target mean to satisfy mean constraint~\eqref{price2}.

Notably, two main challenges exist in~\eqref{eq:grad}:
\begin{enumerate}
    \item The response behavior of ECs should be accurately estimated in~\eqref{eq:map}. 
    However, this task is a complex long-horizon temporal sequence modeling problem, 
    which poses significant challenges in both accuracy and generalizability.
    \item Given $\pi^{\mathrm{EC},k}$, $\partial C_s/\partial \pi^{\mathrm{EC},k}$ has a closed form regarding the scenario value, but
    the gradient computation in~\eqref{eq:grad} requires numerous scenarios to achieve accurate risk quantification,
    which leads to a heavy computational burden.
\end{enumerate}

In the next section, inspired by recent advances in quantum computing,
we propose a quantum learning and estimation approach to address the above two challenges,
thus facilitating efficient DNs-ECs coordinated operation.

\section{Methodology of \\Quantum Learning and Estimation }\label{sec:methodology}
In this section, 
we present a quantum learning and estimation approach.
Specifically, a Q-TCN-LSTM model is proposed to learn the incentive-response characteristics of ECs. 
Subsequently, a quantum estimation method is developed to accelerate the gradient computation across numerous scenarios.

Note that
this study focuses on exploring the theoretical advantages of quantum learning and estimation. Given that current NISQ-era quantum hardware is constrained by noise, decoherence, and limited accessibility, practical deployment of the proposed methods is not yet feasible. Therefore, the subsequent analysis is performed on a simulated ideal quantum device, where external disturbances are excluded and hardware-level improvements (e.g., quantum error correction) are not considered.

\subsection{Quantum Learning}
Generally, QNNs are built on trainable variational quantum circuits (VQCs), 
which consist of layered arrangements of parameterized quantum gates, leveraging 
quantum properties such as superposition and entanglement.
In this part, we first present the fundamental structure of VQCs.
Then, we introduce the proposed Q-TCN-LSTM model.

\subsubsection{Variational Quantum Circuits}
\begin{figure}[htbp]
  \setlength{\abovecaptionskip}{-0.1cm}  
  \setlength{\belowcaptionskip}{-0.1cm}   
  \centering
  \includegraphics[width=1.00\columnwidth]{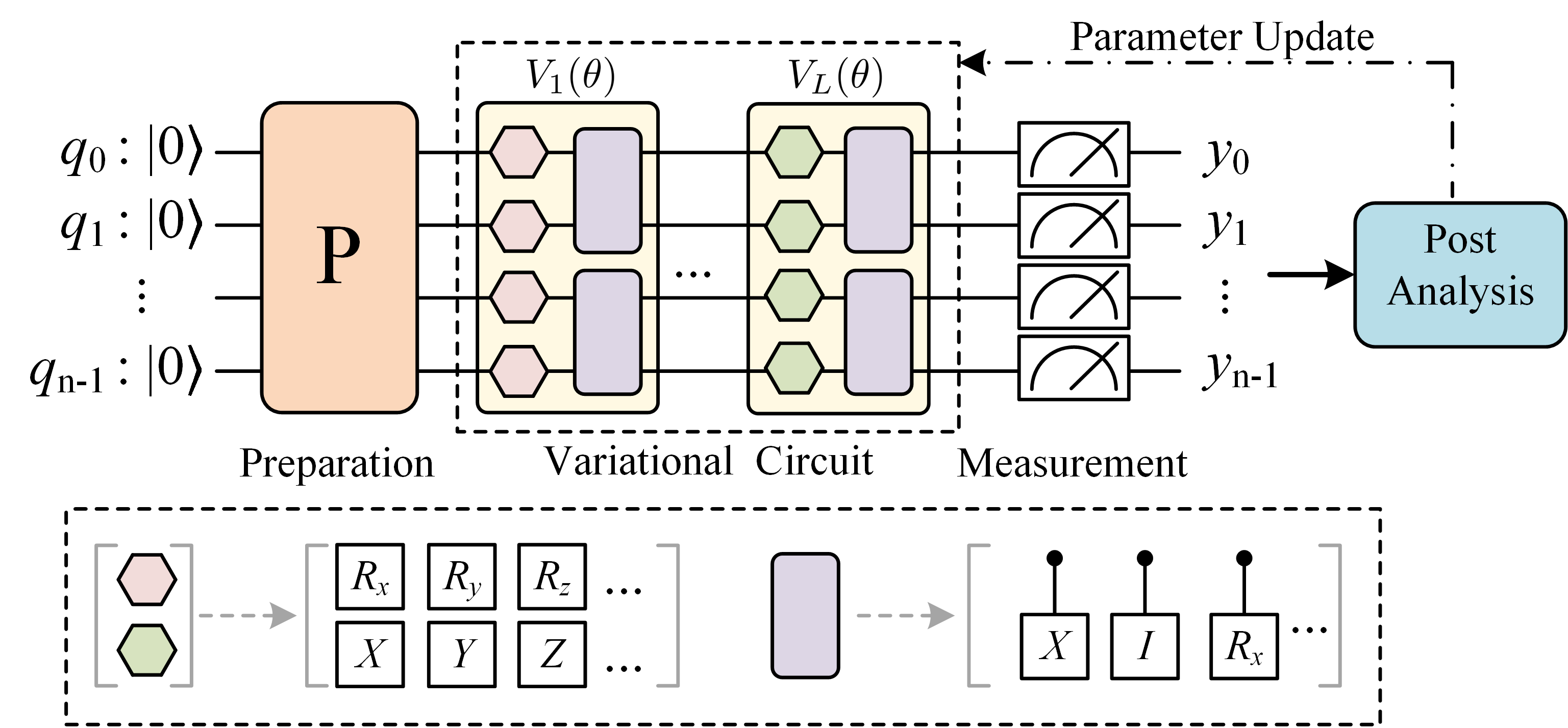}
  \caption{General structure of variational quantum circuits.}
  \label{fig:VQC}
\end{figure}

A VQC is generally composed of three main components: 
quantum state preparation, parameterized variational circuits and quantum measurement, 
as illustrated in Fig.~\ref{fig:VQC}.

Firstly, the quantum state preparation block $\mathbf{P}$ embeds the classical data into quantum states:
\begin{equation}
  \left| \psi_0 \right\rangle_n = \sum_{i=0}^{2^n-1} \alpha_i \left| i \right\rangle_n,
\end{equation}
where $\left| \psi_0 \right\rangle_n$ represents the initial quantum state of $n$ qubits; 
$\left| i \right\rangle_n$ is the $i$-th basis state in the $2^n$-dimensional Hilbert space,
associated with amplitude $\alpha_i$ and a measurement probability of $\alpha_i^2$,
which satisfies $\sum_{i=0}^{2^n-1} \alpha_i^2 = 1$.
At this stage, VQC encodes exponentially many states via quantum superposition, 
offering computational advantages in high-dimensional problems.
Besides, quantum entanglement enables a single quantum operation to act on all $2^n$ superposition states, thus significantly improving computational efficiency.

Secondly, the parameterized variational circuit block $V(\boldsymbol{\theta})$ applies a sequence of complex unitary transformations on $\left| \psi_0 \right\rangle_n$. 
Single-qubit gates (e.g., Pauli-X gate, rotation gate) 
and multi-qubit gates (e.g., controlled-NOT gate, controlled-rotation gate) are used to manipulate the quantum state,
collectively inducing quantum entanglement and further exploiting the superposition across qubits. 
This process enables the circuit to represent intricate nonlinear dependencies inherent in the input data. 
Usually, the variational circuit is composed of $L$ layers, 
each consisting of a set of parameterized or fixed single-qubit gates and entangling gates, forming the overall unitary:
\begin{equation}
  V(\boldsymbol{\theta}) = \prod_{l=1}^{L} (V^{\mathrm{ent}}_{(l)} \cdot V^{\mathrm{rot}}_{(l)}(\boldsymbol{\theta}_{(l)})),
\end{equation}
where $V^{\mathrm{ent}}_{(l)}$ is the entangling gate at layer $l$, 
and $V^{\mathrm{rot}}_{(l)}(\boldsymbol{\theta}_{(l)})=\bigotimes_{i=1}^{n} R_{(x,y,z)}(\theta_{(l),i})$ 
is the rotation gate at layer $l$.
With the variational circuit, the quantum state is transformed as:
\begin{equation}
  \left| \psi(\boldsymbol{\theta}) \right\rangle_n = V(\boldsymbol{\theta}) \left| \psi_0 \right\rangle_n.
\end{equation}

Finally, the measurement block performs projective measurements on $\left| \psi(\boldsymbol{\theta}) \right\rangle_n$:
\begin{equation}\label{q_observe}
  \left\langle \mathcal{O} \right\rangle_{\boldsymbol{\theta}} = \left\langle \psi(\boldsymbol{\theta}) \right| \mathcal{O} \left| \psi(\boldsymbol{\theta}) \right\rangle_n.
\end{equation}
where $\left\langle \psi(\boldsymbol{\theta}) \right|$ is the conjugate transpose of $\left| \psi(\boldsymbol{\theta}) \right\rangle$.
Equation~\eqref{q_observe} denotes the expectation value of the observable
Hermitian operator $\mathcal{O}$, 
resulting in classical scalar values for subsequent analysis.

The measurement output $\left\langle \mathcal{O} \right\rangle_{\boldsymbol{\theta}}$ represents a highly nonlinear transformation of the input, 
implicitly capturing complex correlations encoded by the VQC.

Regarding parameters optimization, classical optimization algorithms (e.g., Adam) are employed to compute the gradients of the predefined loss function (e.g., MSE for regression tasks) with respect to the trainable parameters, where gradients are approximated via methods such as the Parameter-Shift Rule. Parameters are iteratively updated to minimize the loss until convergence.

\subsubsection{Structure of Q-TCN-LSTM}
\begin{figure}[!h] 
  \setlength{\abovecaptionskip}{-0.1cm}  
  \setlength{\belowcaptionskip}{-0.1cm}   
  \centering
  \includegraphics[width=0.8\columnwidth]{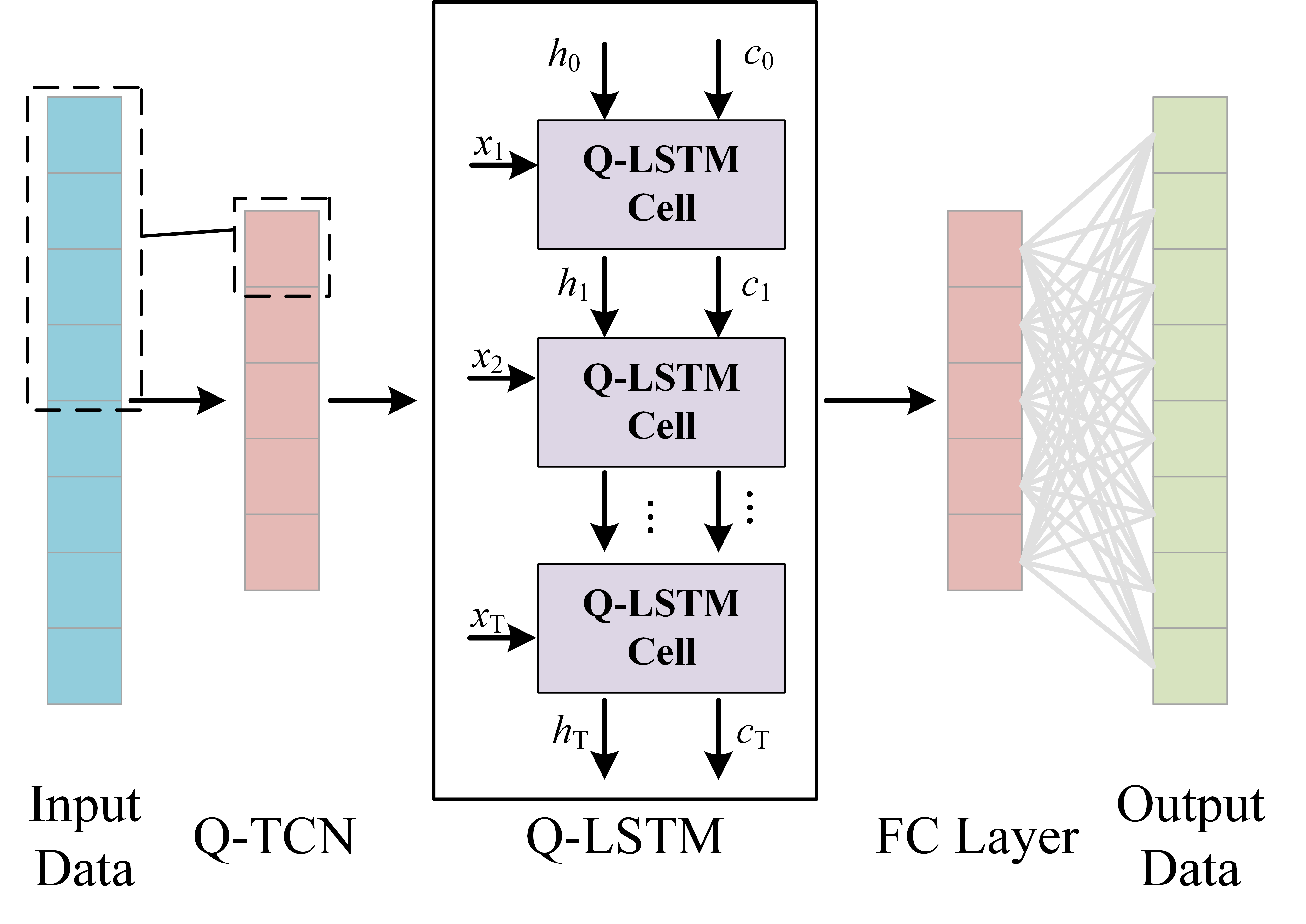}
  \caption{Structure of Q-TCN-LSTM.}
  \label{fig:QTCNLSTM}
\end{figure}

In response to time-varying electricity price signals, ECs exhibit both short-term (spanning several hours) response 
characteristics and long-term (over the scheduling horizon) globally optimal decision-making characteristics. 
Therefore, a short and long term coupling relationship exists between the response of ECs and price incentives.
To capture this coupling, 
we propose a quantum-enhanced hybrid architecture, Q-TCN-LSTM,
which integrates the advantages of 
temporal convolutional networks and long short-term memory networks. 
Specifically, TCNs can efficiently extract short-term features from time series data through causal and dilated convolution mechanisms,
and can adjust the kernel size to control the receptive field.
Then, LSTM can store long-term information in its cell states and selectively remember or forget information through gating mechanisms, 
effectively capturing long-term dependencies in time series.
In the Q-TCN-LSTM model, the convolutional window of Q-TCN slides along the input time series, continuously extracting short-term coupling information within the sequence. Meanwhile, the original long sequence is transformed into a shorter sequence with a more compact dimension, which is then fed as input to the subsequent Q-LSTM layer to enable in-depth mining of data features.
Fully connected (FC) layers then map the extracted features to the final output.
The overall structure of the proposed Q-TCN-LSTM model is shown in Fig.~\ref{fig:QTCNLSTM}.
Notably, the structure of Q-TCN-LSTM (e.g., number of qubits, VQC depth, number of feature extraction layers) can be flexibly adjusted for different tasks.

\subsubsection{Q-TCN for Short-Term Temporal Feature Extraction}
\begin{figure}[htbp]
  \setlength{\abovecaptionskip}{-0.1cm}  
  \setlength{\belowcaptionskip}{-0.1cm}   
  \centering
  \includegraphics[width=1.00\columnwidth]{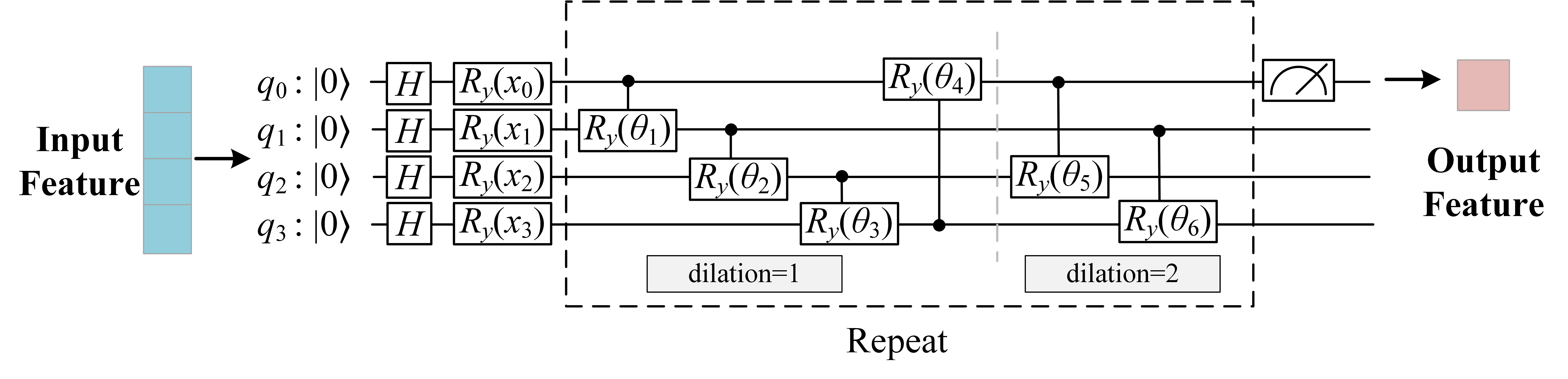}
  \caption{Structure of 4-qubit VQC for Q-TCN.}
  \label{fig:QTCN}
\end{figure}

We integrate VQCs into the causal convolutions to construct a Q-TCN model.
For instance,
a 4-qubit VQC for Q-TCN is illustrated in Fig.~\ref{fig:QTCN}.

Q-TCN utilizes a temporal rotational encoding scheme, 
where each input value $x_t$ is mapped to the angle of a rotation gate. 
This encoding preserves local temporal patterns while mitigating potential information loss 
commonly associated with sparse sampling in classical TCNs.
Then, 
a parameterized variational circuit with skip-connected entanglement gates performs entanglement
across temporally aligned qubits,
constructing temporal dependencies without deep stacking or activations in classical TCNs. 
Finally, quantum measurements produce output features that capture interactions across multiple time steps.
Quantum properties such as superposition and entanglement enhance 
Q-TCN's expressive power with significantly fewer trainable parameters, 
effectively reducing model redundancy.

\subsubsection{Q-LSTM for Long-Term Temporal Feature Extraction}
\begin{figure}[htbp]
  \setlength{\abovecaptionskip}{-0.1cm}  
  \setlength{\belowcaptionskip}{-0.1cm}   
  \centering
  \includegraphics[width=0.9\columnwidth]{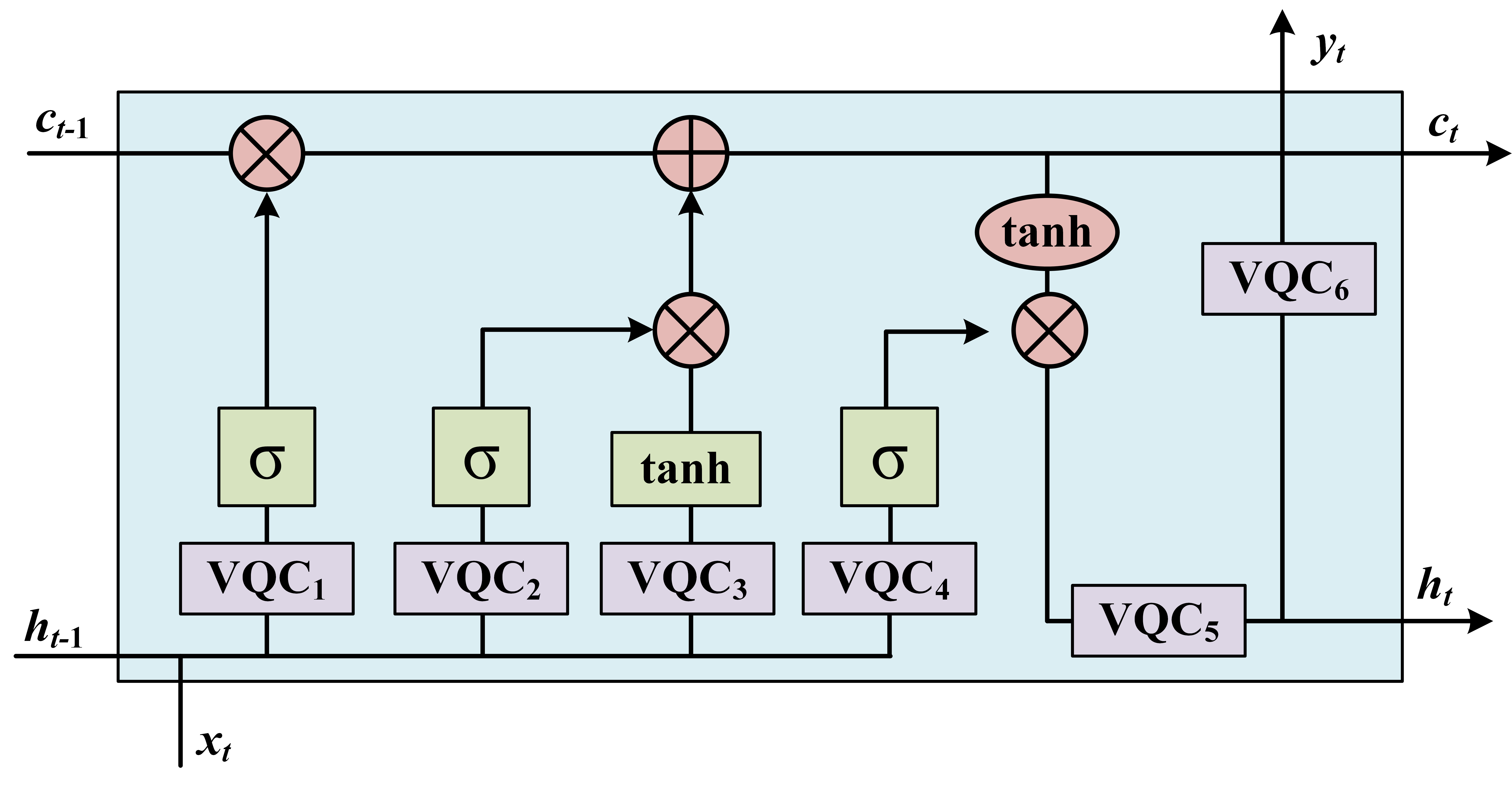}
  \caption{Structure of a Q-LSTM cell.}
  \label{fig:QLSTM}
\end{figure}

Based on the classical LSTM architecture, we construct a Q-LSTM model 
by replacing the gate mechanisms with VQCs,
as illustrated in Fig.~\ref{fig:QLSTM}.

The corresponding forward pass equations in a compact form are as follows:
\begin{subequations}
  \begin{align}
    v_t &= [h_{t-1}, x_t] , \\
    f_t &= \sigma(VQC_1(v_t)),\\
    i_t &= \sigma(VQC_2(v_t)),\ o_t = \sigma(VQC_4(v_t)) , \\
    c_t &= f_t \otimes c_{t-1} + i_t \otimes \tilde{c}_t,\ \tilde{c}_t = \tanh(VQC_3(v_t)) , \\
    h_t &= VQC_5(o_t \otimes \tanh(c_t)),\ y_t = VQC_6(h_t) ,
    \end{align}
\end{subequations}
where 
$\sigma$ and $\tanh$ blocks represent the sigmoid and the hyperbolic tangent activation functions, respectively. 
$x_t/h_t/c_t/y_t$ are the input, hidden state, cell state, and output at time $t$.
$\otimes$ represents element-wise multiplication.
$[\cdot]$ represents vector concatenation.
$f_t$, $i_t$, $o_t$ represent the forget gate, input gate, and output gate in the LSTM model, respectively.

\subsection{Quantum Estimation}
In this part, we first present the general methodology of quantum amplitude estimation.
Then in light of the specific problem of this paper, we provide the detailed application scheme of QAE:

\subsubsection{General Methodology of QAE}
Classical approaches typically rely on MC simulation for uncertainty estimation, 
which has been verified to require $\mathcal{O}(1/\varepsilon^2)$ samples to achieve an estimation error of $\varepsilon$.
This leads to substantial computational overhead.
Alternatively, quantum amplitude estimation,
leveraging quantum entanglement and amplitude amplification,
can theoretically achieve the same precision level $\varepsilon$ with only $\mathcal{O}(1/\varepsilon)$ queries~\cite{QAE} under certain conditions\footnote{
Such conditions include the efficient preparation of the initial quantum state, small hardware noise in real devices, and stringent precision requirements, among others.},
offering a theoretical quadratic speedup compared to MC.

The QAE procedure generally involves three steps:

\textbf{Step 1-Quantum State Preparation:} 
A distribution $\mathcal{X}$ is divided into $N=2^n$ intervals,
and is embedded into a quantum state:
\begin{equation}\label{q_state_prepare}
  |\psi\rangle_n = \sum_{i=0}^{N-1} \sqrt{p_i} |i\rangle_n,
\end{equation}
where $|i\rangle_n$ represents a basis state of $n$ qubits,
$p_i$ is the probability of state $|i\rangle_n$ with $\sum_{i=0}^{N-1} p_i = 1$.
The essence of preparing a quantum state embedded with distribution $\mathcal{X}$ lies in mapping probabilities to a quantum superposition state via amplitude encoding.
In the uniform case where $p_i = 1/N\  \forall i\in [0,\dots,N-1]$,~\eqref{q_state_prepare} can be efficiently prepared using Hadamard gates.
In more complex cases, it is often necessary to design sophisticated quantum state preparation circuits: for instance, constructing continuous Gaussian distributions using quantum arithmetic circuits~\cite{QGaussian1,QGaussian2}, or approximating complex distributions through parameterized quantum circuits~\cite{QGaussian3}.

\textbf{Step 2-Quantum State Rotation:}
Assuming a function $f:\{0,1,...,2^n-1\}\rightarrow [0,1]$,
a controlled rotation operator $\mathcal{F}$ is constructed such that 
\begin{equation}\label{eq:rotation}
  \mathcal{F}|i\rangle_n|0\rangle = \sqrt{1-f(i)}|i\rangle_n |0\rangle + \sqrt{f(i)}|i\rangle_n |1\rangle.
\end{equation}

Applying $\mathcal{F}$ to $|\psi\rangle_n$ yields
\begin{equation}\label{eq:QAE}
  \mathcal{F}|\psi\rangle_n|0\rangle = \sum_{i=0}^{N-1} \sqrt{(1-f(i))p_i} |i\rangle_n|0\rangle + \sum_{i=0}^{N-1} \sqrt{f(i)p_i} |i\rangle_n|1\rangle.
\end{equation}

In the implementation of $\mathcal{F}$, 
the auxiliary qubit undergoes a rotation $R_y(\theta)$ with $\theta = 2\arcsin(\sqrt{f(i)})$,
which means that $\theta$ is a nonlinear function of state $i$.
Depending on the trade-off between circuit complexity and estimation accuracy, we propose two potential implementations of $\mathcal{F}$.

\begin{figure}[htbp] 
  \setlength{\abovecaptionskip}{-0.1cm}  
  \setlength{\belowcaptionskip}{-0.1cm}   
  \centering
  \includegraphics[width=1.00\columnwidth]{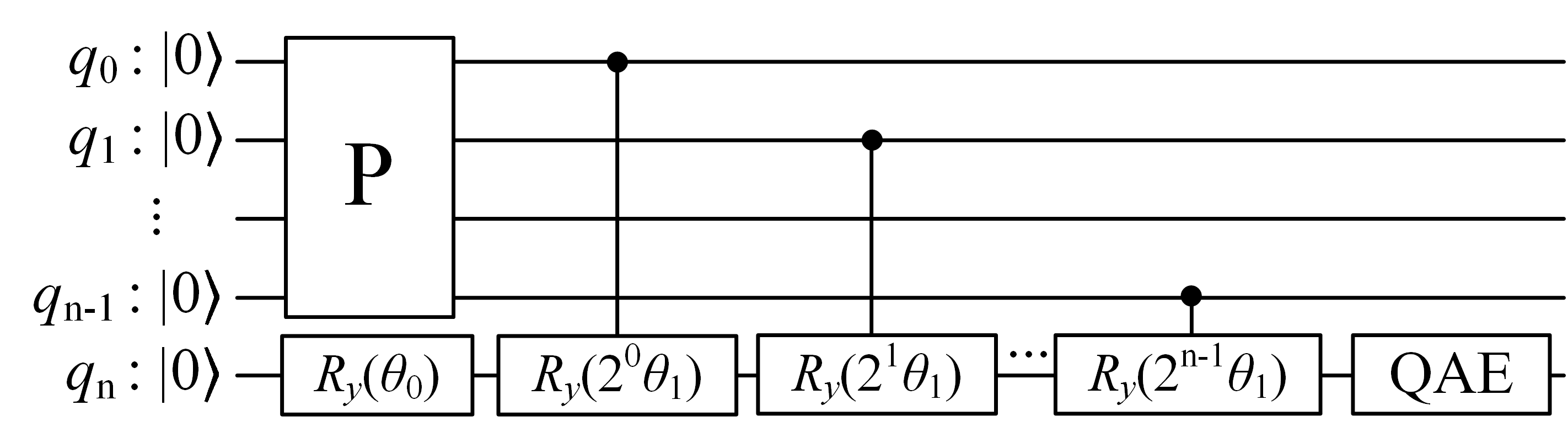}
  \caption{The approximated implementation of $\mathcal{F}$, referred to as Circuit-1.}
  \label{fig:QAEappro}
\end{figure}

Circuit-1 employs analytical approximations of the nonlinear rotation, as illustrated in Fig.~\ref{fig:QAEappro}.
The core idea is to find a linear formulation with respect to $i$ that efficiently approximates $\arcsin(\sqrt{f(i)})$ such that $\theta \approx 2(ai + b)$. For example, when $x$ is sufficiently small, we have $\arcsin(\sqrt{x}) \approx ax+b$, where $a$ and $b$ are determined using the least-squares method.
Under this condition, the controlled rotation operator $\mathcal{F}$ can be implemented using a sequence of parameterized controlled-RY gates, with low implementation complexity.
This approximation benefits from easy implementation but 
may introduce certain approximation errors.

\begin{figure}[htbp] 
  \setlength{\abovecaptionskip}{-0.1cm}  
  \setlength{\belowcaptionskip}{-0.1cm}
  \centering
  \includegraphics[width=1.00\columnwidth]{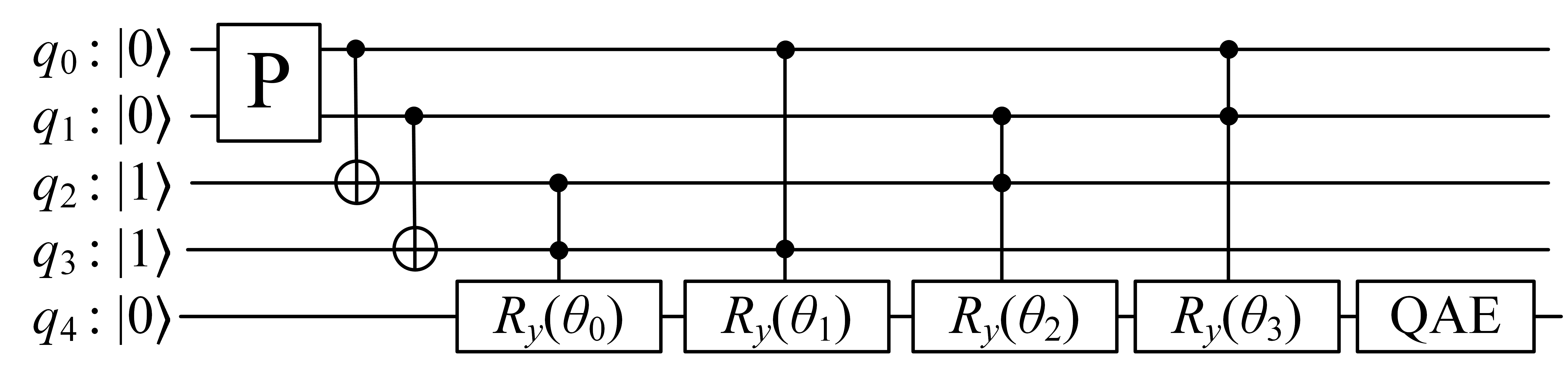}
  \caption{The precise implementation of $\mathcal{F}$, referred to as Circuit-2.}
  \label{fig:QAE}
\end{figure}

Circuit-2 computes the exact nonlinear rotation angle corresponding to each quantum basis state.
To this end, 
we design a rotation-angle-adaptive circuit architecture. 
In this scheme, 
$n$ original qubits are initialized to $|0\rangle_n$, 
$n$ reversal qubits to $|1\rangle_n$, 
and an auxiliary qubit to $|0\rangle$. 
The reversal qubits are kept in bitwise negated states of the original qubits by applying CNOT gates.
Subsequently, 
multi-controlled rotation gates are employed to apply precise $R_y(\theta)$ operations to the auxiliary qubit, 
conditioned on the combined states of both the original and reversal qubits. 
For example, 
in the case of $n=2$, 
the proposed adaptive rotation circuit is illustrated in Fig.~\ref{fig:QAE}. 
When the original qubits are in $\left| 10 \right\rangle$, 
the reversal qubits are automatically in $\left| 01 \right\rangle$ due to the CNOT configuration. 
A controlled rotation gate is then activated based on this unique qubit configuration to apply the corresponding rotation $R_y(\theta_2)$ to the auxiliary qubit.
This implementation achieves high-precision realization of the controlled rotation operator $\mathcal{F}$. 
However, 
this accuracy comes at the cost of increased circuit depth and significant gate complexity, 
as $2n+1$ qubits are required and 
$2^n$ controlled rotation gates are added to the circuit,
which may hinder scalability for larger $n$.

\textbf{Step 3-Quantum Amplitude Estimation:}
Based on~\eqref{eq:QAE}, we can apply the QAE algorithm to estimate the amplitude of the ancillary qubit being in $|1\rangle$,
corresponding to the expectation of $\sum_{i=0}^{N-1} f(i)p_i$.
Commonly adopted QAE algorithms include canonical QAE, iterative QAE and maximum likelihood QAE. The percentage estimation error relative to the ground truth is defined as $\epsilon=$ $|estimation-truth|/truth \times  100\%$.

For brevity, we summarize the fundamental idea of QAE~\cite{QAE} as follows.
Consider a quantum system with $n$ qubits initialized in $|0\rangle_n$, 
and an ancillary qubit initialized in the state $|0\rangle$. 
Let $\mathcal{A}$ be a quantum unitary operator that:
\begin{equation}\label{eq:amplitude}
    \mathcal{A}|0\rangle_n|0\rangle = \sqrt{1-a}|\psi_0\rangle_n|0\rangle + \sqrt{a}|\psi_1\rangle_n|1\rangle.
\end{equation}
where $a\in [0,1]$ is unknown.  $|\psi_0\rangle_n$ and $|\psi_1\rangle_n$ are two normalized states.
QAE aims to estimate $a$ with high probability.
To achieve this goal, QAE employs the Grover operator 
$\mathcal{Q} = \mathcal{A}\mathcal{S}_0\mathcal{A}^{\dagger }\mathcal{S}_{\mathcal{X}}$,  
where $\mathcal{S}_0$ is the operator for selective phase inversion of the initial state, 
and $\mathcal{S}_{\mathcal{X}}$ is the operator for phase marking of the ``good'' state.
The original amplitude estimation problem can be mapped to a phase estimation problem by iterative action on $\mathcal{Q}$. 
The amplitude $a\in [0,1]$ is encoded as a rotation angle $\theta \in [0,\frac{\pi}{2}]$,
satisfying $a = \sin^2(\theta)$, and renders~\eqref{eq:amplitude} as:
\begin{equation}
  \mathcal{A}|0\rangle_n|0\rangle = \cos(\theta)|\psi_0\rangle_n|0\rangle + \sin(\theta)|\psi_1\rangle_n|1\rangle.
\end{equation}

By performing $k$ iterations of the Grover operator $\mathcal{Q}$,
the quantum state is transformed to:
\begin{equation}
  \mathcal{Q}^k\mathcal{A}|0\rangle_n|0\rangle = \cos((2k+1)\theta)|\psi_0\rangle_n|0\rangle + \sin((2k+1)\theta)|\psi_1\rangle_n|1\rangle.
\end{equation}

Thus, the probability of measuring the ancillary qubit in state $|1\rangle$ is $\sin^2((2k+1)\theta)$.
When $k$ is chosen to be $\lfloor \frac{\pi}{4\theta} \rfloor $, we have $\sin^2((2k+1)\theta) \approx 1$.
Here, we achieve the amplitude amplification of the target state $|1\rangle$.
It has been verified that the estimated value $\tilde{a}$ satisfies
$|\tilde{a} - a| \leq O(\frac{1}{M^{\mathrm{Q}}})$
with probability of at least $\frac{8}{\pi^2}$. 
$M^{\mathrm{Q}}$ is the number of oracle queries.
QAE algorithm has been verified to be able to achieve a quadratic speedup compared to the $O(\frac{1}{\sqrt{M^{\mathrm{S}}}})$ 
convergence rate of classical MC methods~\cite{MC1,QAE} under relatively ideal conditions, 
where $M^{\mathrm{S}}$ is the number of MC simulation samples.

\subsubsection{Method of applying QAE to (4)}
Next, we propose a quantum estimation method for efficient computation of~\eqref{eq:grad}.
Based on the LinDistFlow model in~\eqref{DNs_opt}, variables such as nodal voltage and line power flow exhibit a linear relationship with the uncertain variables (i.e., RES generation, load demand) and $P^{\mathrm{EX,EC}}$~\cite{lindistflow,lindistflow1}.
Thus, gradient terms $\partial ( \pi^{\mathrm{EC}}_{t}  P^{\mathrm{EX,EC}}_{i,t} )/\partial \pi^{\mathrm{EC}}_t$ and $\partial   P^{\mathrm{EX,DN}}_{s,t} /\partial \pi^{\mathrm{EC}}_t$
can be easily computed.
Regarding the penalty terms with $[\cdot]_+$, note that the distributions of $V_{i,s,t}$ and $P_{ij,s,t}$ can be accurately determined based on the known distribution of uncertain variables,
the relationship between these penalty gradient terms (e.g., $\frac{\partial[V_{i,s,t} - \overline{V}_{i}]_+}{\partial \pi^{\mathrm{EC}}_t}$) and scenario index $i$ presents a piecewise characteristic.
Note that the computational procedures for these penalty terms are essentially the same.
Take one penalty term  
$\sum_{s=1}^{S} \gamma_s \frac{\partial[V_{i,s,t} - \overline{V}_{i}]_+}{\partial \pi^{\mathrm{EC}}_t}$ as an example, the core steps are as follows:
a) encoding the distribution of $V_{i,t}$ using a quantum state preparation
circuit, 
b) adopting specific controlled rotation circuits to model the penalty term,
where Circuit-1 in Fig.~\ref{fig:QAEappro} can be utilized for implementation
if the linear approximation error is acceptable,
otherwise Circuit-2 in Fig.~\ref{fig:QAE} can be employed if the original $f(i)$ is retained.
Such a choice can be determined based on the trade-off between accuracy and computational complexity. 
c) applying a QAE block to the auxiliary qubit to obtain the estimation result.
Similarly, other penalty terms can also be computed efficiently. 
Based on the above analysis,~\eqref{eq:grad} can be computed efficiently by applying QAE.

\section{Case Studies}\label{sec:case}

\subsection{Settings}
A modified IEEE 33-bus distribution network integrated with DGs, uncertain loads, and multiple ECs
is adopted as the test system.
One WT and two PVs are connected to buses 5, 11, and 31, respectively.
Three uncertain loads are located at buses 10, 12, and 14, 
while the remaining loads are assumed to be deterministic.
An industrial EC, a commercial EC, and a residential EC are connected to buses 7, 15, and 18, respectively. 
The voltage magnitude is restricted as $|V_i| \in [0.90,1.10]\ \text{(p.u.)},\ \forall i \in \varOmega^\mathrm{N}$. 
$\Delta t$ = 15 min and $T$ = 96.
$\lambda$ = 1 and $\alpha$ = 0.95.
The VQC for Q-TCN model uses 4 qubits and 3 layers. 
The VQC in the Q-LSTM cell uses a hidden size of 3 and 2 recurrent layers.
The learning rate is set to 0.01. 
Three datasets comprising 10,240 incentive-response samples from three ECs are 
utilized to train the proposed mapping model.
QNNs are implemented using PennyLane, Jax, and Flaxnn libraries.
QAE is implemented using Qiskit library.
All computations are implemented in Python 3.10.11,
conducted on a Windows 11 64-bit operating system with an Intel Core i9-13900HX @ 2.30 GHz processor and 16 GB RAM.

\subsection{Performance of Quantum Learning}
\begin{figure}[htbp]
  \setlength{\abovecaptionskip}{-0.1cm}  
  \setlength{\belowcaptionskip}{-0.1cm}
  \centering
  \includegraphics[width=1.00\columnwidth]{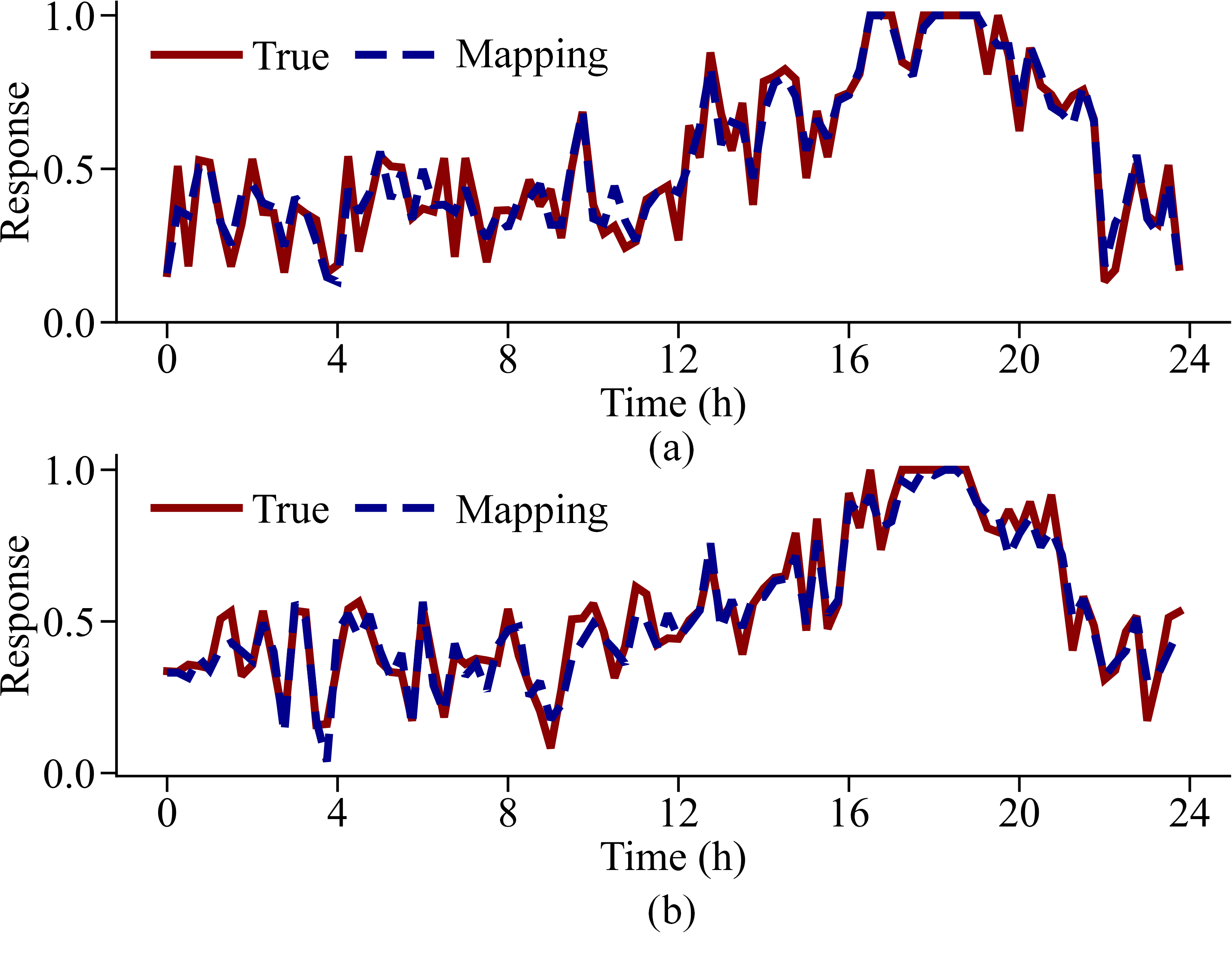}
  \caption{Mapping results for two commercial EC samples.}
  \label{fig:fit_result}
\end{figure}

We evaluate the performance of Q-TCN-LSTM in terms of accuracy and computational efficiency. 
In terms of convergence speed,
compared to the classical TCN-LSTM, which converges in about 25 epochs with an MSE of 0.013, 
Q-TCN-LSTM converges faster within 10 epochs and achieves a lower MSE of 0.004. 
These results demonstrate the improved learning ability of Q-TCN-LSTM,
benefiting from quantum advantage.
Furthermore,
we visualize two randomly selected response trajectories of the commercial EC from the test set, as shown in Fig.~\ref{fig:fit_result}. 
The mapping results (blue dotted line) closely match the true values (red solid line) across both high-value and low-value ranges
with minimal mapping deviation, demonstrating the model's high mapping precision and temporal generalization performance.

\begin{figure}[htbp]
  \setlength{\abovecaptionskip}{-0.1cm}  
  \setlength{\belowcaptionskip}{-0.1cm}
  \centering
  \includegraphics[width=1.00\columnwidth]{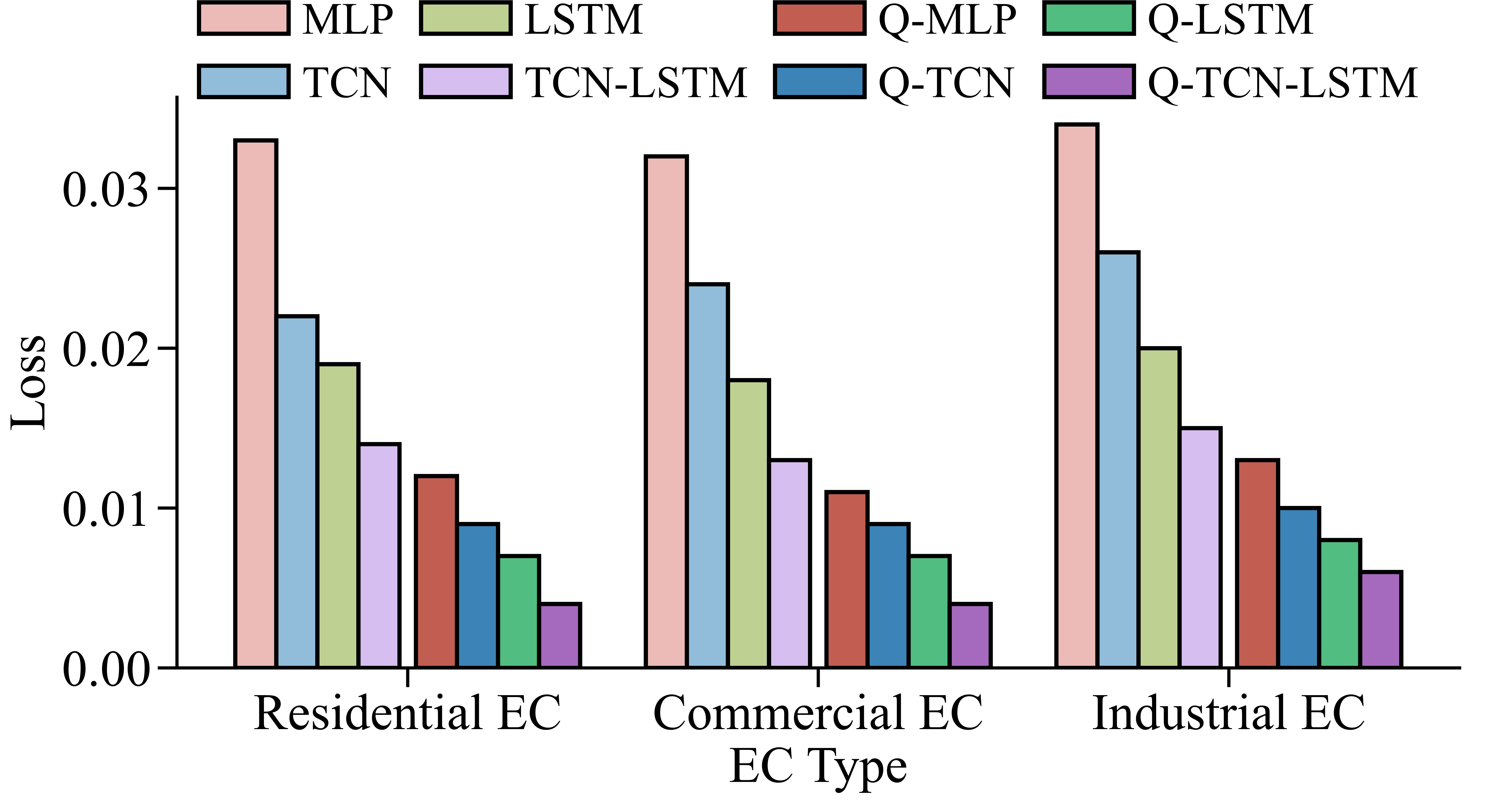}
  \caption{Performance comparison of multiple models.}
  \label{fig:loss_compare}
\end{figure}

In addition, 
a comprehensive comparative experiment is conducted to further evaluate the superiority of Q-TCN-LSTM. 
Classical neural networks (Multi-Layer Perceptron (MLP), TCN, LSTM, TCN-LSTM) 
and quantum neural networks (Q-MLP, Q-TCN, Q-LSTM) are included for comparison,
as shown in Fig.~\ref{fig:loss_compare}.

As shown in Fig.~\ref{fig:loss_compare},
QNNs consistently outperform their classical counterparts across all three ECs, 
achieving lower MSE and demonstrating superior function approximation capabilities. 
This advantage can be attributed to the unique properties of quantum computing. 
Specifically, quantum entanglement provides QNNs with powerful nonlinear expressive capabilities, 
enabling the effective modeling of complex relationships between pricing incentives and ECs' responses.

Among the QNN architectures evaluated, 
Q-TCN-LSTM achieves the highest mapping accuracy. 
This is primarily due to its hybrid architecture, 
which integrates the strengths of both Q-TCN and Q-LSTM networks. 
In Q-TCN, the integration of quantum convolutional layers enables the model 
to efficiently capture fine-grained temporal dependencies within short time sequences. 
The Q-LSTM leverages quantum superposition to significantly increase its memory capacity. 
This allows the model to maintain strong correlations between distant data points, leading to a better understanding of long-range context.
By leveraging the strengths of both Q-TCN and Q-LSTM, 
the Q-TCN-LSTM architecture is capable of capturing both short-term and long-term temporal dependencies inherent in the EC response patterns, 
thereby enhancing its overall performance.

Moreover, we compare the computational efficiency of Q-TCN-LSTM with other models.
Due to the current limitations in quantum hardware availability, 
all quantum operations are simulated on classical processors.
Thus, for the following time efficiency computation,
we separately measure two types of time indicators:
\textit{simulation runtime} ($\tau^{\mathrm{s}}$), which is the running time of task execution based on a quantum simulator on a classical computer; 
and \textit{quantum runtime}  ($\tau^{\mathrm{q}}$) which is the estimated time for an ideal quantum hardware device to run the corresponding task.
Specifically, we adopt the formula in~\cite{quantumtime} to approximate the quantum runtime:
$T^\mathrm{Q} \approx T^\mathrm{P} + T^\mathrm{G} \cdot D + T^\mathrm{M}$,
where $T^\mathrm{P}$ denotes state preparation time,
$T^\mathrm{G}$ represents the average duration of a single quantum gate,
$D$ is the maximal circuit depth,
and $T^\mathrm{M}$ is the quantum measurement time. 
Empirical values from~\cite{quantumtime} suggest that $T^\mathrm{P}+T^\mathrm{M}=1\  \mu s$ and $T^\mathrm{G}=10\  ns$,
which indicates that the time for state preparation and measurement accounts for a significant portion of the total time
\footnote{
Note that the runtime of quantum algorithms on actual hardware is affected by multiple factors, such as the specific qubit platform, measurement methods, and error correction mechanisms, among others. Specific quantum runtime measurements are also reported in~\cite{RN942,fan2025calibrating, ni2025swap}.
}.
Thus, one execution of the convolution operation in the VQC in Q-TCN takes around 0.03 $s$ of simulation runtime and 1.1 $\mu s$ of quantum runtime,
which is much less than that of classical convolution and activation operation (around 0.017 $s$).
This is because quantum parallelism enables performing multiple operations simultaneously, 
significantly boosting computational speed. 
In terms of model size,
the TCN-LSTM model has a total of 37353 parameters, while the Q-TCN-LSTM model uses only 94 parameters, 
approximately only 0.25\% of the TCN-LSTM model's parameter count.
Q-TCN-LSTM outperforms the classical TCN-LSTM while using far fewer parameters and operating at a much higher speed, 
demonstrating its high computational efficiency.

To investigate the impacts of noise on model performance, which is inherent in practical quantum hardware in the NISQ era,
we simulated multiple noise levels during the model training process (theoretically [0, 1], with practical simulations focusing on [0, 0.2]). 
The relationship between the mapping performance (measured by MSE) and the noise level is:
\begin{table}[htbp]
\centering
\caption{Model Performance under Different Noise Levels}
\begin{tabular}{lccccc}
\toprule
Noise Level & 0 & 0.05 & 0.1 & 0.15 & 0.2 \\
\midrule
MSE     & 0.004 & 0.005 & 0.006 & 0.01 & 0.12 \\
\bottomrule
\end{tabular}
\end{table}

This result indicates that with the increase in noise level, the MSE shows a gradual upward trend. 
The core reason is that noise in quantum devices interferes with the stability of quantum states, 
leading to errors in the quantum computing process and thus reducing the model mapping accuracy. 
Therefore, in practical quantum hardware, quantum error correction strategies are needed to suppress noise interference and ensure computing performance.

\subsection{Performance of Quantum Estimation}
\begin{figure}[htbp] 
  \setlength{\abovecaptionskip}{-0.1cm}  
  \setlength{\belowcaptionskip}{-0.1cm}
  \centering
  \includegraphics[width=1.00\columnwidth]{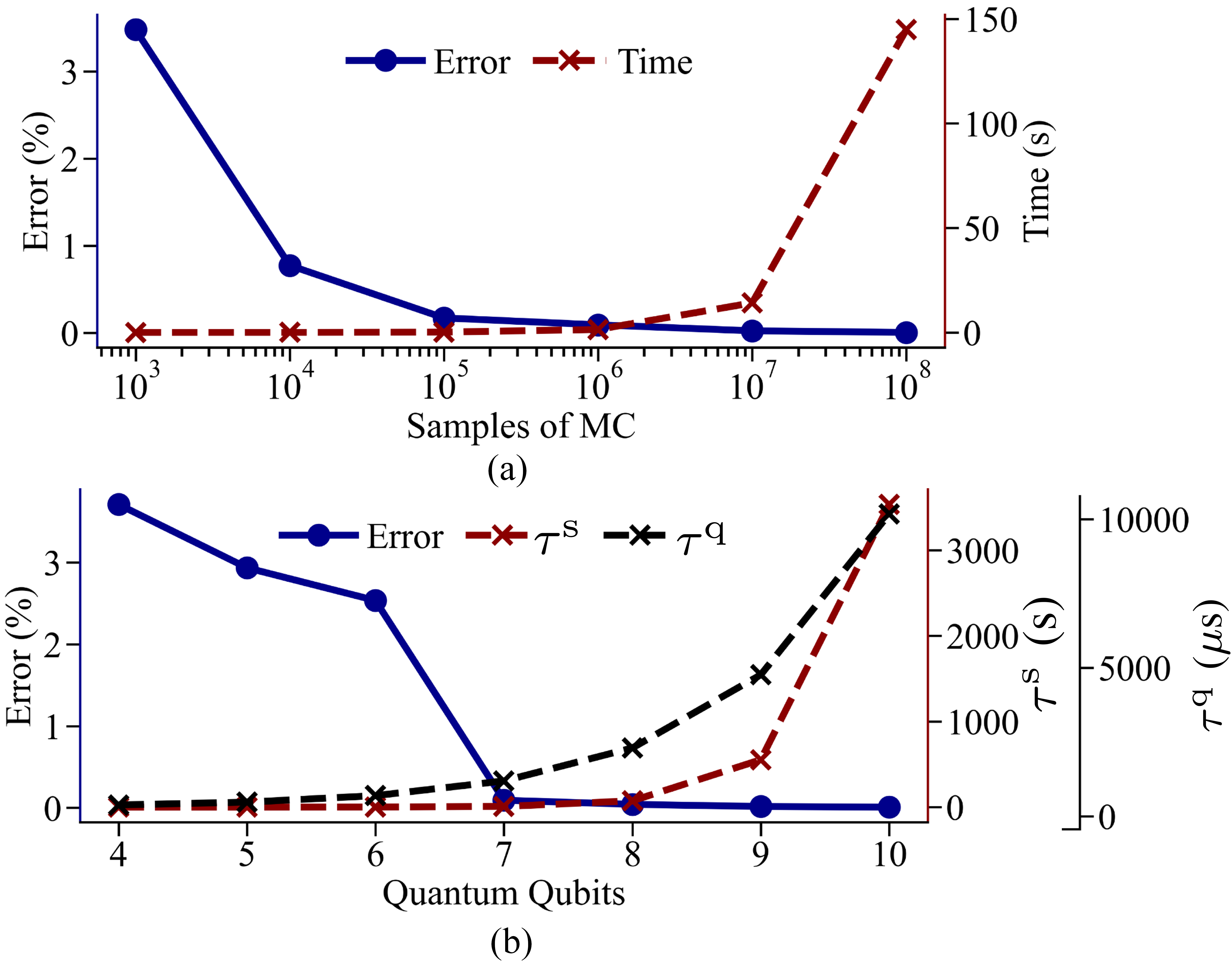}
  \caption{Comparison of MC and QAE: (a) MC, (b) QAE with Circuit-2.}
  \label{fig:QAEMC}
\end{figure}

To validate the computational efficiency of quantum estimation methods, 
we compare the estimation error (expressed as a percentage) and the corresponding computation time including both simulation runtime and quantum runtime
under varying numbers of MC samples and QAE qubits. 
The detailed results are illustrated in Fig.~\ref{fig:QAEMC}.

As shown in Fig.~\ref{fig:QAEMC}, 
MC achieves a low estimation error of approximately $0.1\%$ with a substantial sample size of $10^6$,
while QAE attains a comparable error level of $0.09\%$ with only 7 qubits.
This demonstrates a substantial reduction in required computational resources.
Moreover, QAE demonstrates a marked advantage in time efficiency. 
Under a similar error threshold (MC with $10^6$ samples vs. QAE with 7 qubits), 
MC requires approximately 1.41 s of computation time, 
while QAE completes the task in approximately 
8.74 s of simulation runtime and around 920 $\mu s$ of quantum runtime.
This efficiency gain is primarily attributed to the quadratic speedup in query complexity offered by QAE.
Furthermore, given the exponential growth of quantum circuit execution capacity with increasing qubits, 
the computational advantage of QAE is expected to become even more pronounced in large-scale estimation problems. 
Overall, these findings underscore the superiority of QAE in terms of both estimation accuracy and computational efficiency.

In addition, 
the performance of QAE using the approximate circuit (Circuit-1)
and the precise circuit (Circuit-2) is compared in terms of estimation error and computational time,
as shown in Table~\ref{tab:QAE_comparison}.
The results indicate that Circuit-1 exhibits a notable increase in estimation error compared to Circuit-2.
This is primarily due to the linear approximation error. 
In contrast, Circuit-2 can represent such functions with higher precision but at the cost of increased complexity:
it requires $2n+1$ qubits and $2^{n}$ quantum controlled rotation gates, 
resulting in substantially higher execution time.
Therefore, in practical applications, 
the choice between the two circuit designs should be guided by the functional characteristics of the estimation target 
and the required trade-off between accuracy and computational complexity. 
\begin{table}[htbp]
  \setlength{\abovecaptionskip}{-0.1cm}  
  \setlength{\belowcaptionskip}{-0.1cm}  
    \setlength{\tabcolsep}{0.166cm} 
    \centering
    \caption{Performance Comparison of Circuit-1 and Circuit-2}
    \label{tab:QAE_comparison}
    \begin{tabular}{cccccccc}
        \toprule
        \multirow{2}{*}{Circuit}&\multirow{2}{*}{Metric} & \multicolumn{6}{c}{$n$} \\ \cmidrule(lr){3-8}
        & & 5 & 6 & 7 & 8 & 9 & 10  \\ \midrule
        \multirow{3}{*}{Circuit-1} 
        & $\epsilon$ (\%)    & 9.123 & 6.688 & 5.459 & 4.834 & 4.534 & 4.379  \\  
        & $\tau^{\mathrm{s}}$ (s) & 0.088 & 0.210 & 0.617 & 1.607 & 3.876 & 9.741  \\
        & $\tau^{\mathrm{q}}$ ($\mu$s) & 5.64 & 6.52 & 7.44 & 8.36 & 9.28 & 10.21  \\ 
        \multirow{3}{*}{Circuit-2} 
        & $\epsilon$ (\%)     & 2.933 & 2.534 & 0.093 & 0.042 & 0.017 & 0.008 \\ 
        & $\tau^{\mathrm{s}}$ (s) & 0.225 & 1.714 & 8.742 & 73.53 & 553.7 & 3534.4  \\
        & $\tau^{\mathrm{q}}$ ($\mu$s) & 176 & 403 & 920 & 2073 & 4636 & 10271  \\ 
        \bottomrule
    \end{tabular}
\end{table}

\subsection{Result of DNs-ECs Coordinated Operation}

\begin{figure}[htbp]
  \setlength{\abovecaptionskip}{-0.1cm}  
  \setlength{\belowcaptionskip}{-0.1cm}   
  \centering
  \includegraphics[width=1.00\columnwidth]{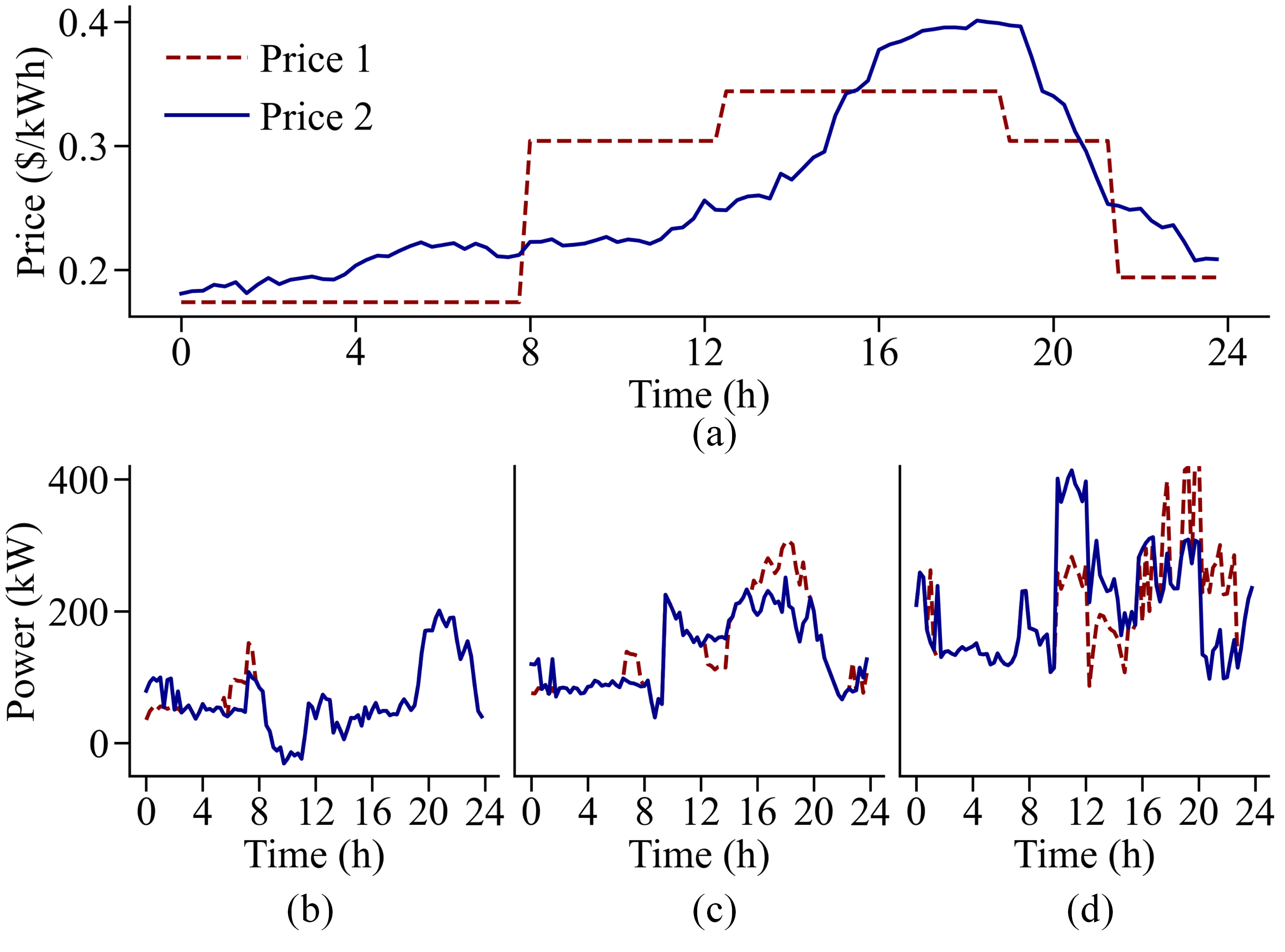}
  \caption{Results of incentive price and responses of ECs: (a) Electricity trading price, (b) Residential EC, (c) Commercial EC, (d) Industrial EC.}
  \label{fig:ECDR}
\end{figure}

In the optimization problem of DNs,
leveraging the accurate incentive-response mapping of ECs enabled by Q-TCN-LSTM,
along with the accurate gradient estimation enabled by QAE, 
gradient descent can be efficiently employed.
The bilevel DNs-ECs coordinated operation optimization results are shown in Fig.~\ref{fig:ECDR}. 

As illustrated in Fig.~\ref{fig:ECDR}, 
the optimized pricing scheme (Price 2) exhibits more time-varying characteristics compared to the baseline fixed pricing (Price 1). 
This enables better alignment with the comprehensive operation patterns of multiple ECs.
Fig.~\ref{fig:ECDR} indicates that Price 2, compared with Price 1,
exhibits a relatively higher price during hours 16-20, while presenting a relatively lower price during hours 8-16. This pricing pattern closely matches the underlying net load dynamics of the ECs: 
during hours 8-16, high penetration of renewable generation leads to lower net loads, 
making it suitable to reduce electricity prices to encourage load shifting. 
Conversely, during hours 16-20,
the renewable generation decreases, 
and the load demand of ECs slightly increases,
and thus the net load demand reaches a relative peak demand interval, 
where elevated prices are employed to suppress load aggregation and incentivize shifting the load from this period to other periods.
The effectiveness of this pricing scheme is further reflected in user response behavior, 
as shown in Fig.~\ref{fig:ECDR}(b)-(d).
Following price optimization, all three ECs exhibit varying degrees of net load shifting, 
particularly during the high-price hours 16-20, resulting in a noticeable reduction in peak net loads. 
Among them, industrial users demonstrate the most substantial response due to their large and adjustable load capacity, 
indicating high sensitivity to price signals. 
Commercial users show moderate responsiveness, while residential users exhibit limited flexibility due to the inherent rigidity 
of their consumption patterns and fewer controllable loads. 
This coordinated response behavior effectively reduces the negative impact of ECs on the DNs.
After optimizing the pricing strategy,
the average penalty cost of voltage violation decreases from 50.9 to 21.8,
indicating the capacity of user-side response to alleviate potential risks and enhance distribution network reliability.
Overall, leveraging the advantages of quantum computing,
the DNs-ECs coordinated operation optimization problem can be solved efficiently.

\section{Conclusion}\label{sec:conclusion}
In this paper, a quantum learning and estimation approach is proposed to address the bilevel DNs-ECs coordinated operation problem.
By leveraging quantum properties including quantum superposition and entanglement,
the proposed Q-TCN-LSTM model can efficiently learn the incentive-response characteristics of ECs,
realizing accurate estimation of user electricity consumption behavior under limited information.
Additionally, 
QAE is introduced to facilitate evaluations with numerous scenarios, 
achieving a great speedup over classical MC methods,
significantly reducing computational resource consumption and execution time. 
In this way, the proposed quantum learning and estimation approach 
efficiently facilitates the coordinated operation between DNs and ECs to generate better pricing strategies to regulate the net load profiles of ECs,
thereby reducing the operational risk introduced to the DNs.
As illustrated by the case studies,
the proposed quantum-enhanced approach demonstrates superior performance in terms of both mapping accuracy and computational efficiency.

Although quantum learning and estimation methods have demonstrated certain theoretical advantages over their classical counterparts, the practical application of quantum computing in power systems is constrained by current quantum hardware limitations.  
First, hardware performance and scalability are inadequate. Existing quantum processors provide only a few hundred qubits with short coherence times and high noise sensitivity, limiting their ability to handle large-scale power system tasks. Although error correction may improve reliability, it requires significant additional resources.  
Second, the quantum–classical interface restricts efficiency. Classical data must be encoded into quantum states and later decoded through measurement, causing potential information loss and being limited by qubit availability.  
Third, quantum hardware is costly and difficult to integrate with existing infrastructures. Devices require extreme operating conditions and are incompatible with current classical data systems, such as SCADA, without substantial modifications.  
Hence, advancing practical quantum hardware is critical to converting theoretical advantages into engineering value, enabling accurate risk assessment and optimized operation in power systems.

\bibliographystyle{IEEEtran}
\bibliography{IEEEabrv,Quantum}

\begin{IEEEbiography}[{\includegraphics[width=1in,height=1.25in,clip,keepaspectratio]{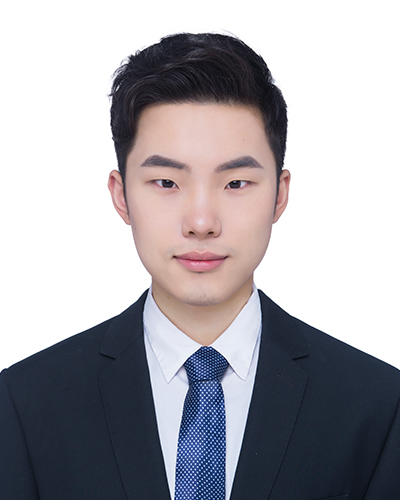}}]{Yingrui Zhuang}
    (Graduate Student Member, IEEE) was born in 1999. He received the B.S. degree in electrical engineering in 2021 from Tsinghua University,
    Beijing, China, where he is currently working toward the Ph.D. degree in electrical engineering with the
    State Key Laboratory of Power System Operation and Control, Department of Electrical Engineering.
    His research interests include risk analysis of distribution systems with high renewable energy resources penetration 
    and data-driven distributed energy resources coordination and optimization.
\end{IEEEbiography}

\begin{IEEEbiography}[{\includegraphics[width=1in,height=1.25in,clip,keepaspectratio]{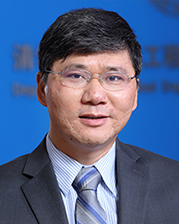}}]{Lin Cheng}
    (Senior Member, IEEE) was born in 1973. He received the B.S. degree in electrical engineering from Tianjin University, Tianjin, China, in 1996 
    and the Ph.D. degree from Tsinghua University, Beijing, China, in 2001. 
    He is currently a tenured Professor with the Department of Electrical Engineering, Tsinghua University, where he is also the Deputy
    Director of the State Key Laboratory of Power System Operation and Control.
    His research interests include operational reliability evaluation and application of power systems, operation optimization of distribution
    systems with flexible resources, and perception and control of uncertainty in wide-area measurement systems.
    He is also a Fellow of IET.
\end{IEEEbiography}

\begin{IEEEbiography}[{\includegraphics[width=1in,height=1.25in,clip,keepaspectratio]{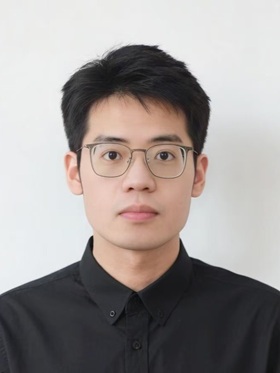}}]{Yuji Cao} received the B.E. and M.S. degrees from The Chinese University of Hong Kong, Shenzhen, China, in 2021 and 2023, respectively. He is currently pursuing the Ph.D. degree in the Department of Mechanical and Automation Engineering, The Chinese University of Hong Kong, Hong Kong, China. His research focuses on privacy-preserving algorithms with applications in smart grids.
\end{IEEEbiography}

\begin{IEEEbiography}[{\includegraphics[width=1in,height=1.25in,clip,keepaspectratio]{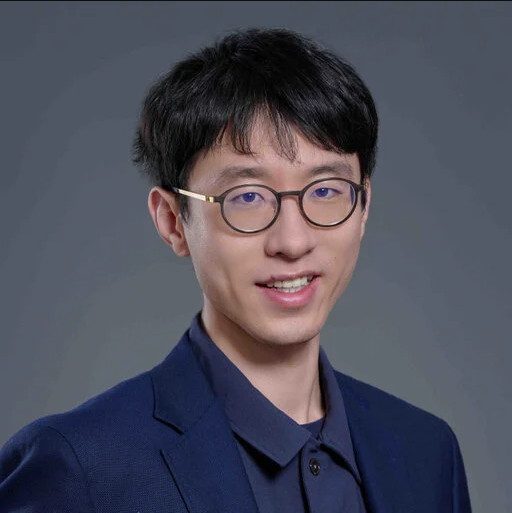}}]{Tongxin Li}
(Member, IEEE) received the
Ph.D. degree from the Department of Computing
and Mathematical Sciences, Caltech, Pasadena, CA,
USA in 2022.
He is an Assistant Professor and a Presidential
Young Fellow with the School of Data Science
(SDS), The Chinese University of Hong Kong,
Shenzhen, China. Prior to joining SDS, in 2022,
he interned as an Applied Scientist at Amazon
Web Services (AWS), Boston, MA, USA, in 2020
and 2021. His current research interests include
interdisciplinary topics in power systems, control, machine learning, and
online algorithms.
Dr. Li was a recipient of the 2022 ACM SIG Energy Doctoral Dissertation
Award Honorable Mention.
\end{IEEEbiography}

\begin{IEEEbiography}[{\includegraphics[width=1in,height=1.25in,clip,keepaspectratio]{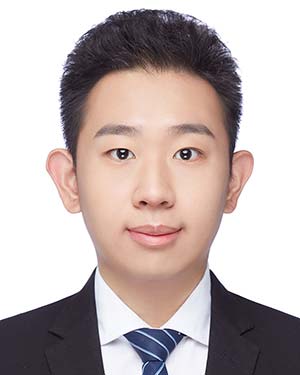}}]{Ning Qi}
    (Member, IEEE) was born in 1996. He received the B.S. degree in electrical engineering from Tianjin University, Tianjin, China, in 2018 
    and the Ph.D. degree in electrical engineering from Tsinghua University, Beijing, China, in 2023.
    He was a Visiting Scholar with the Technical University of Denmark, Lyngby, Denmark, in 2022.
    He was a Research Associate in Electrical Engineering with Tsinghua University, in 2024. He is currently a Postdoctoral Research
    Scientist in earth and environmental engineering with Columbia University, New York, NY, USA.
    His research interests include behavior modeling, optimization under uncertainty, and market design for power systems 
    with generalized energy storage. 
    He is also the Youth Editorial Board Member of Power System Protection and Control, and the Guest Editor of Processes and Frontiers in Energy Research.
\end{IEEEbiography}

\begin{IEEEbiography}[{\includegraphics[width=1in,height=1.25in,clip,keepaspectratio]{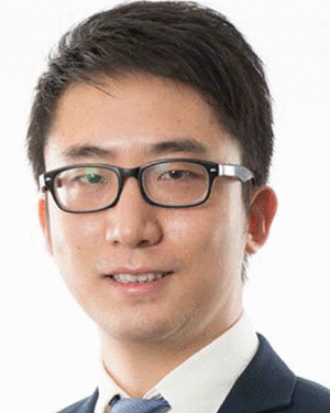}}]{Yan Xu} (Senior Member, IEEE) received the B.E.
and M.E. degrees from the South China University
of Technology, Guangzhou, China, in 2008 and
2011, respectively, and the Ph.D. degree from the
University of Newcastle, Callaghan, NSW, Australia,
in 2013. He joined the Nanyang Technological
University (NTU), Singapore, with the Nanyang
Assistant Professorship after Postdoc training. He
was promoted to an Associate Professor in early
2021 and appointed as the Cham Tao Soon Professor
in engineering (an endowed professorship named
after NTU’s founding president) in early 2024. He is currently the Director
with the Center for Power Engineering, NTU, and the Co-Director of
Singapore Power Group-NTU Joint Lab. His research interests include power
system stability and control, microgrids, and data-analytics for smart grid
applications. He was the recipient of the University of Sydney Postdoctoral
Fellowship for postdoctoral training. He was an Associate Editor for the
IEEE TRANSACTIONS ON SMART GRID and the IEEE TRANSACTIONS
ON POWER SYSTEMS, the Chairman of IEEE Power and Energy Society
Singapore Chapter in 2021 and 2022, respectively, and a General Co-Chair
of the 11th IEEE ISGT-Asia Conference in 2022.
\end{IEEEbiography}

\begin{IEEEbiography}[{\includegraphics[width=1in,height=1.25in,clip,keepaspectratio]{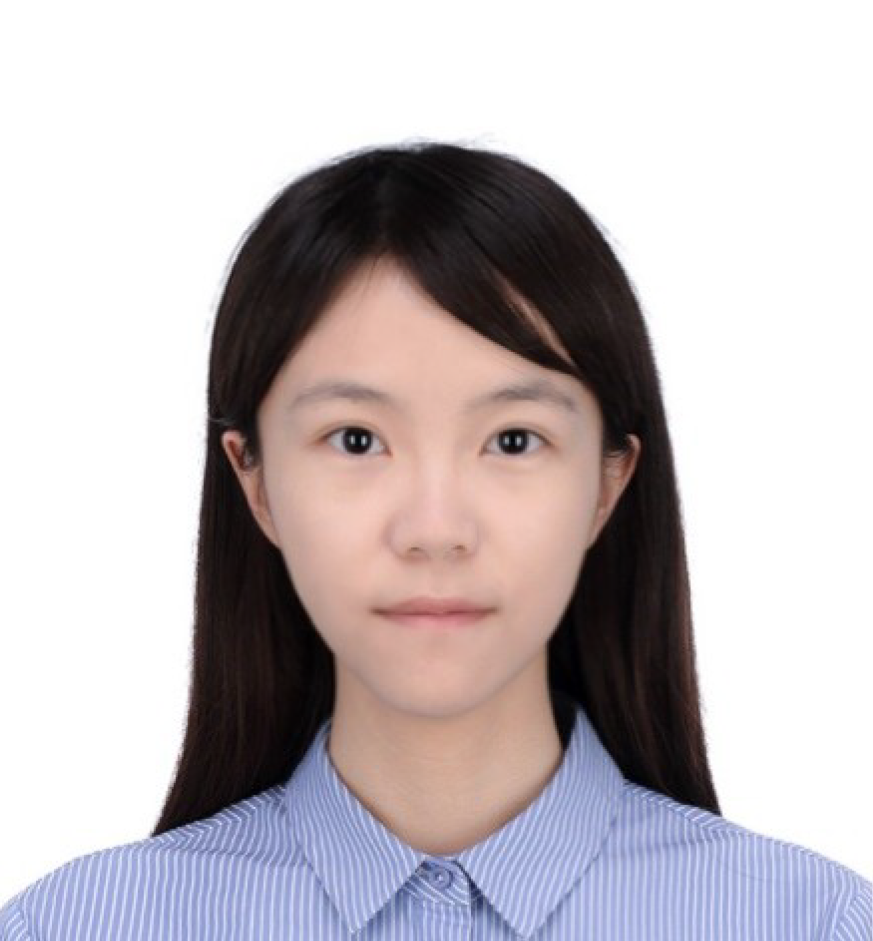}}]{Yue Chen}
 (Senior Member, IEEE) received the
B.E. degree in electrical engineering from Tsinghua
University, Beijing, China, in 2015, the B.S. degree in economics from Peking University, Beijing,
China, in 2017, and the Ph.D. degree in electrical engineering from Tsinghua University, in 2020.
She was a visiting student with California Institute
of Technology, Pasadena, CA, USA from 2018 to
2019. She is currently a Vice-Chancellor Assistant
Professor with the Department of Mechanical and
Automation Engineering, The Chinese University
of Hong Kong, Hong Kong, China. Her research interests include robust
optimization, game theory, and trustworthy AI with application to smart grids.
She is currently an Associate Editor for IEEE TRANSACTIONS ON SMART
GRID, IEEE TRANSACTIONS ON SUSTAINABLE ENERGY, and IEEE POWER
ENGINEERING LETTERS.
\end{IEEEbiography}

\end{document}